\documentclass{article}
\usepackage[utf8]{inputenc}
\usepackage[a4paper,margin=2.5cm,bmargin=2.5cm]{geometry}

\usepackage{graphicx}
\usepackage{subfig}
\graphicspath{{./figures/}}
\usepackage{dcolumn}
\usepackage{bm}
\usepackage[mathlines]{lineno}
\usepackage[export]{adjustbox}
\usepackage{color}
\usepackage{tabularx,array,booktabs}
\usepackage{caption}
\usepackage{mathtools}
\usepackage{amsmath,amssymb}
\usepackage{soul}
\usepackage{lipsum}
\usepackage{varioref}
\usepackage{algorithm,algpseudocode}
\usepackage{ulem}
\usepackage{siunitx}
\usepackage{comment}
\usepackage{upgreek}
\usepackage[dvipsnames]{xcolor}
\definecolor{ColoreSeba}{HTML}{FF8300}

\renewcommand{\phi}{\varphi}
\newcommand{\bb}[1]{\mathbf{#1}}
\newcommand{\bs}[1]{\boldsymbol{#1}}
\newcommand{\R}{\mathbb{R}}
\newcommand{\C}{\mathbb{C}}
\newcommand{\LL}{\mathcal{L}}
\newcommand\norm[1]{\left\lVert#1\right\rVert}
\DeclarePairedDelimiter{\norma}{\lVert}{\rVert}

\newcommand{\Comsol}{COMSOL Multiphysics\textsuperscript{\textregistered}}

\begin{document}

\title{Design of arbitrarily shaped acoustic cloaks through PDE-constrained optimization satisfying sonic-metamaterial design requirements}

\date{}

  \author{Sebastiano Cominelli \thanks{Sebastiano Cominelli is a PhD Candidate at Politecnico di Milano, Department of Mechanical Engineering, Milano 20133, Italy (e-mail: sebastiano.cominelli@mail.polimi.it), Corresponding author.} \and
  Davide E. Quadrelli \thanks{Davide E. Quadrelli is a PhD Candidate at Politecnico di Milano, Department of Mechanical Engineering, Milano 20133, Italy (e-mail: davidee.quadrelli@polimi.it).} \and
    Carlo Sinigaglia\thanks{Carlo Sinigaglia is a PhD Candidate at Politecnico di Milano, Department of Mechanical Engineering, Milano 20133, Italy (e-mail: carlo.sinigaglia@polimi.it)} \and
    Francesco Braghin\thanks{Prof. Francesco Braghin is Full Professor at Politecnico di Milano, Department of Mechanical Engineering, Milano 20133, Italy (e-mail: francesco.braghin@polimi.it). }
  }

\maketitle
\begin{abstract}
We develop an optimization framework for the design of acoustic cloaks, with the aim of overcoming the limitations of usual transformation-based cloaks in terms of microstructure complexity and shape arbitrarity of the obstacle. 
This is achieved by recasting the acoustic cloaking design as a nonlinear optimal control problem constrained by a linear elliptic partial differential equation.
In this setting, isotropic material properties' distributions realizing the cloak take the form of control functions and a system of first-order optimality conditions is derived accordingly.  Such isotropic media can then be obtained in practice with simple hexagonal lattices of inclusions in water.

For this reason, the optimization problem is directly formulated to take into account suitable partitions of the control domain Two types of inclusions are considered, and long-wavelength homogenization is used to define the feasible set of material properties that is employed as a constraint in the optimization problem.
In this manner, we link the stage of material properties optimization with that of microstructure design, aiming at finding the optimal implementable solution. As a test benchmark, cloaking of the silhouette of a ship is considered, for various frequencies and directions of incidence. The resulting cloak is numerically tested via coupled structural/acoustic simulations.

\end{abstract}

\section{Introduction}
\label{intro}
Inspired by the development of Transformation Theory, the quest for the implementation of invisibility devices has spread during the last decade over diverse research fields \cite{kadic2015experiments}, in which governing partial differential equations have been shown to be invariant under coordinate transformations. Started in electromagnetism \cite{pendry2006controlling, leonhardt2006optical}, this theory has indeed unlocked the possibility to achieve perfect concealment from detection in acoustics \cite{cummer2007one,chen2007acoustic, norris2008acoustic}, elastodynamics
\cite{norris2011elastic, norris2012hyperelastic}, surface water waves \cite{farhat2008broadband}, heat conduction \cite{schittny2013experiments} and even matter waves \cite{zhang2008cloaking}. The beauty and power of this analytical method stands in the fact that the obtained cloak is theoretically exact for all frequencies and incoming directions of the probing incident radiation.

In acoustics \cite{norris2008acoustic, norris2009acoustic}, it has also been shown that the solution of the problem is not unique in terms of material parameters distributions: \textit{inertial cloaks} \cite{torrent2008acoustic, pendry2008acoustic, popa2011experimental, zigoneanu2014three} are made with anisotropic inertial properties, pure \textit{pentamode cloaks} \cite{chen2015latticed, chen2017broadband} are obtained with solids exhibiting singular anisotropic elasticity tensors, while the most general acoustic cloak can comprise both mass and elasticity anisotropy. On the flip side, however, this material distributions are hard to achieve in practice, and one is often forced to resort to complex microstructures designed by homogenization-based optimization techniques \cite{layman2013highly, kadic2012practicability}, to obtain the required anisotropic material behavior. More than that, analytical solutions are available for simple geometries only, such as the axisymmetric case \cite{gokhale2012special}, and the literature dealing with arbitrarily shaped cloaks based on transformation theory is limited and almost entirely restricted to the inertial cloak case \cite{li2012homogeneous, li2018non, li2019two,chen2016design, quadrelli2021}. 

Several attempts have been made to overcome such restrictions and allow for simplified design, for instance using quasi-conformal cloaks \cite{li2008hiding}, in which the transformation is specifically constructed in such a way that anisotropy is avoided in the obtained material distributions. However, the geometries that allow for application of this technique are limited and the cloak should in principle comprise the overall space, thus a truncation is required that makes the solution not exact.

\textit{Scattering cancellation}, instead, is an alternative technique that relies on surrounding the target with a distribution of small obstacles, in such a way that the resulting multiple scattering solution has no influence on the incident field. Such distribution of scatterers can be obtained by setting \textit{a priori} their number and shape and optimizing for their location either with evolutionary algorithms \cite{garcia2011acoustic} or with gradient based optimization \cite{amirkulova2017acoustic}. Increased degrees of freedom can instead be considered in the optimization if not only the location, but the shape also is not fully determined a priori: in this case one can use parametric optimization of Bezier shapes \cite{lu2018acoustic}, or even topology optimization \cite{andkjaer2013topology} which has recently allowed to consider acoustic-elastic interactions \cite{fujii2021acoustic} in the optimal design of acoustic cloaks. The simplicity of construction unlocked by this techniques has also allowed for the design and validation of three dimensional cloaks of axisymmetric obstacles \cite{sanchis2013three}, whose practical demonstration is still lacking when considering classic transformation theory based cloaks. The downside of these methods is that they are inherently narrowband and work for a limited set of incident angles. The broadbandness and the number of working directions can be increased by augmenting the set of cases considered in the cost function, accepting a trade-off between performance and number of working frequencies/directions.

In the search for simplified configurations, one can progressively rely more on optimization and less on model based intuition by exploiting neural networks to compute the physical properties of a set of layers of isotropic homogeneous fluids \cite{ahmed2021deterministic}.

In this paper, we follow another route and reformulate the design phase such that the properties of the cloak are obtained as the solution of a Partial Differential Equation (PDE)-constrained optimization problem, that is an Optimal Control Problem (OCP). The control functions are infinite-dimensional and space-varying fields of material properties that nullify the scattered wave.
The state equation is represented by the inhomogeneous Helmholtz equation \cite{bergmann1946wave} describing the scattered wave in the domain. A similar PDE-constrained optimization framework is considered in \cite{chen2021optimal} with the additional complexity of adding uncertainty in the problem formulation. Instead of considering the wave propagation velocity as control function as in \cite{chen2021optimal}, we consider as separate control functions both the density and bulk modulus fields, and introduce constraints in the optimization for such controls, thus taking into account for the fact that in practical implementations these two parameters can hardly be chosen independently. This in turns allows to derive an elegant and concise expression for the reduced gradient of the cost functional with respect to these two control variables.
More than that, our formulation is intended to facilitate the link between the design of the \textit{macrostructure}, i.e.\ the material property distribution, with that of the \textit{microstructure} that implements via long-wavelength homogenization the required density and bulk modulus, thus unlocking the marriage between the two stages of the design of such two-scale optimization problem. Indeed, the standard approach to implement inhomogeneous material property distributions in acoustic cloaking is to discretize them and fill each resulting sub-domain with an appropriately optimized microstructure \cite{chen2015latticed,chen2017broadband,quadrelli2021experimental}.
This approach leads to sub-optimal solutions depending on the chosen discretization: provided that the sub-domains are sufficiently small compared to the wavelength considered, the wave "feels" a gradient of refraction index that might be different from the required one.
In this work we instead make use of appropriate control basis functions that allow to obtain optimal solutions taking into account the size and shape of the cloak sub-domains at the level of the optimization problem.
Finally, considering inhomogeneous but isotropic material distributions considerably reduces the complexity of the required microstructure, which can be simply obtained considering hexagonal lattices of solid inclusions in the hosting water medium. The manuscript is organized as follows: in the next two sections the optimization problem is introduced, and the optimality conditions are derived. The OCP is then discretized with the Finite Element Method (FEM) in order to allow for numerical solutions and the solution of the usual axisymmetric cloak is shown. In the fourth section, an in deep analysis of the reachable set of homogenized material properties is conducted on simple hexagonal lattices of solid inclusions in water, in order to build a set of constraints for the OCP that allow for practical implementations. In Section 5, such constraints are introduced in the formulation of the problem, and constrained solutions are compared to those obtained previously with the unconstrained problem. Before drawing conclusions, Section 6 deals with the numerical validation of the cloak implemented with the microstructures analyzed in Section 4. A boat-shaped target is considered as an additional case study to validate the method against arbitrarily shaped targets.

\section{Problem Statement}
\begin{figure}
	\centering
	\includegraphics[width=0.5\textwidth]{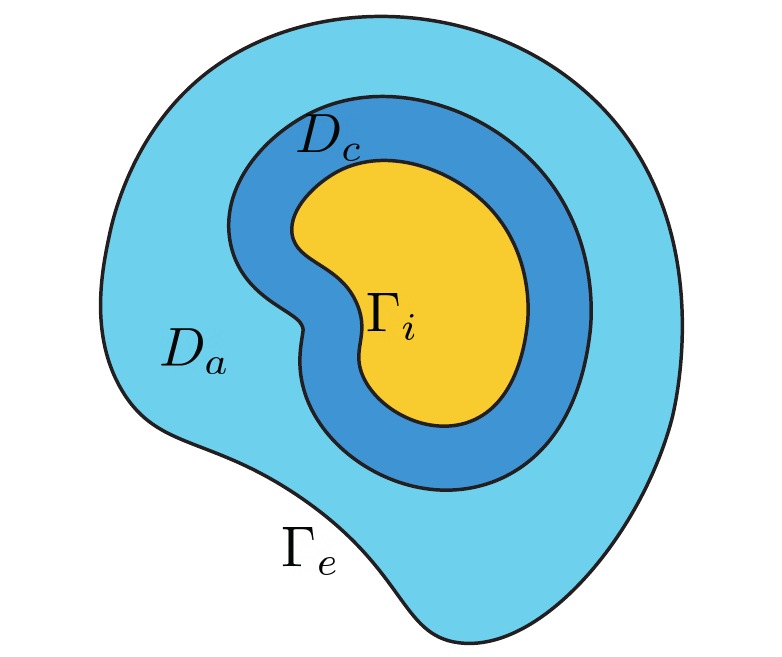}
	\caption[Domain definition]{Schematic representation of the computational domain $\Omega$ with bounday $\Gamma= \Gamma_i \cup \Gamma_e$. The background fluid occupies the subdomain $D_a$, while the cloak occupies the region $D_c$. }
	\label{fig:domain definition}
\end{figure}

We consider a two-dimensional acoustic scattering problem in an inhomogenous medium consisting of water as background fluid and of a cloaking region modeled as an inhomogeneous yet isotropic equivalent fluid. The computational domain $\Omega \subset \R^2$ is divided in two subdomains: $D_c$ is the domain occupied by the cloak, $D_a$ corresponds to the surrounding ambient and it is occupied by the fluid. The domain's boundary is $\Gamma = \partial\Omega = \Gamma_i \cup \Gamma_e$, where $\Gamma_i$ is the obstacle's shape and $\Gamma_e$ the external boundary. The interface between cloak and fluid domains is denoted as $\Gamma_c$ whereas the external boundary $\Gamma_e$ is needed for computational purposes and its role will be detailed in the following. The domain $D_a$ is filled with water with standard physical properties ($\rho_0 = \SI{998}{[\kilo\gram\per\cubic\meter]}$, $\kappa_0=\SI{2.2}{[\mega\pascal]}$). We denote as $\rho_0$ and $\kappa_0$ its properties in the background domain $D_a$. On the other hand, the physical properties in the domain $D_c$ are assumed as control functions and denoted as $\rho$ and $\kappa$. These are considered as function of the space variable $\mathbf{x} \in D_c$. This layout is shown in Figure~\ref{fig:domain definition}. 

When the system is forced by time harmonic waves, the steady-state acoustic pressure $P(\bb{x},t)$ can be separated as $\Re(p(\bb{x})e^{ j \omega t})$, where $\Re\{\cdot\}$ denotes the real component of its argument. The complex amplitude $p(\bb{x})$ satisfies the Helmholtz equation for inhomogeneous media \cite{bergmann1946wave}:
\begin{equation}	\label{eq:Helmholtz frequency}
	\nabla\cdot\Big(a(\bb{x}) \nabla p(\bb{x})\Big) = -b(\bb{x})\omega^2 p(\bb{x}),
\end{equation}
where $p(\bb{x}) \in \C$ is the pressure field phasor, $\omega$ the circular frequency of the forcing wave. The coefficients $a$ and $b$ are defined as $a := \rho^{-1}$ and $b:=\kappa^{-1}$, where $\rho$ is the local mass density and $\kappa$ the local bulk modulus. The definition of $a$ and $b$ will turn out to be useful in manipulating Equation (\ref{eq:Helmholtz frequency}) and setting up the resulting Optimal Control Problem (OCP). \\
The total pressure field can be decomposed into an incident and a scattered field, that is:
\begin{equation}	\label{eq: p_inc + p_sc}
	p(\bb{x})=p_{s}(\bb{x})+p_{i}(\bb{x}).
\end{equation}
where $p_{i}$ is the solution of the Helmholtz equation obtained considering a homogeneous fluid without obstacles. That is $p_{i}$ satisfies Equation \eqref{eq:Helmholtz frequency} with homogeneous properties
\begin{equation}	\label{eq:Helmholtz frequency homogeneous}
	a_0 \Delta p_{i}(\bb{x}) = - b_0 \omega^2  p_{i}(\bb{x}).
\end{equation}
A plane wave solution to Equation~\eqref{eq:Helmholtz frequency homogeneous} is $p_{i}=e^{-j  k_0\,\mathbf{a}\cdot\mathbf{x}}$  where $k_0 = \frac{\omega}{c_{0}}$ is the homogeneous wave number, $c_0=\sqrt{\frac{\kappa_0}{\rho_0}}$ is the undisturbed sound velocity and $\mathbf{a} \in \R^2$ is the unit vector associated to the direction of the incident wave.
Equation~\eqref{eq:Helmholtz frequency} can be rewritten in terms of the scattered pressure $p_s$:
\begin{equation}
\label{eq:state}
	-\nabla \cdot \Big(a(\bb{x})\nabla p_{s}(\bb{x})\Big) - \omega^2 b(\bb{x}) p_{s}(\bb{x}) = \omega^2\Big(b(\bb{x})-b_0\Big)p_{i}(\bb{x}) +\nabla\cdot\Big[\Big(a(\bb{x})-a_0\Big)\nabla p_{i}(\bb{x})\Big]
\end{equation}
which is obtained plugging Equation~\eqref{eq:Helmholtz frequency homogeneous} into Equation~\eqref{eq:Helmholtz frequency} and rearranging the terms. We remark that the incident wave $p_i$ is a datum of the problem. We consider the scattering problem from an infinitely rigid obstacle at the boundary $\Gamma_i$. This in turns specifies a zero normal velocity of the total pressure $p$ as boundary condition, that is:
\begin{align}\label{eq:rigid boundary}
	\nabla p_{s}(\bb{x})\cdot\mathbf{n}(\bb{x}) &= -\nabla p_{i}(\bb{x})\cdot\mathbf{n}(\bb{x}) & \text{on } \Gamma_i,
\end{align}
where $\mathbf{n}$ is the outgoing normal. Note that Equation \eqref{eq:rigid boundary} is a standard inhomogeneous Neumann boundary condition since $p_i$ is completely known.
In order to approximate computationally an unbounded domain we need to guarantee that the scattered wave is outgoing by satisfying the Sommerfeld radiation condition \cite{schot1992eighty}:

\begin{equation}
\label{somm}
\lim_{r \to \infty} \sqrt{r} \Big( \frac{\partial p_s(\bb{x})}{\partial r} + j k_0 p_s(\bb{x}) \Big) = 0
\end{equation}
where $r = \norma{\bb{x}}$. For the sake of simplicity, we substitute Equation~\eqref{somm} with the first-order Bayliss and Turkel approximation for 2D domains, that is \cite{bayliss1980radiation}:
\begin{align}\label{eq:Sommerfeld 1st}
	\nabla p_{s}(\bb{x})\cdot\mathbf{n}(\bb{x}) + \Big(jk_0 + \frac{1}{2R}\Big) p_{s}(\bb{x}) &= 0	& \text{on } \Gamma_e,
\end{align}
where $R$ is the radius of $\Gamma_e$. Note that Equation~\eqref{eq:Sommerfeld 1st} is a homogeneous Robin boundary condition. This approximation guarantees reliable results without increasing the problem complexity (see e.g.\ \cite{SHIRRON1998121}).

\noindent In the following we will omit the explicit dependence on the space variable when it is clear from the context.

\section{The Optimal Control Problem}
\label{sec:OCP}

In this section, the acoustic cloaking problem is formulated as an OCP where the state dynamics consists of the scattered field $p_{s}(\bb{x})$ that solves the linear elliptic PDE \eqref{eq:state}. Space modulated density and bulk modulus in the cloaking region take the role of control functions. Hence, the overall OCP is nonlinear due to the way the control affects the state. 
The cloaking objective is achieved if the intensity of the scattered wave vanishes, that is equivalent to minimize the quadratic objective $\norma{p_{s}(\bb{x})}^2 = \bar{p}_{s}(\bb{x}) p_{s}(\bb{x})$ in the ambient domain $D_a$, where $\bar{p}_{s}(\bb{x})$ represents the complex conjugate of $p_{s}(\bb{x})$. This objective can be encoded in a quadratic cost functional which aims at finding the optimal trade-off minimizing the scattered wave with control functions which deviate as little as possible from the background properties of water. Then, the OCP can be written as follows: 
\begin{align}
	\label{eq:cost functional}
	\min_{v,u,p_{s}} J(v,u,p_{s}) &=  \frac{\lambda_v}{2}\int_{D_c} v^2 \,d\Omega + \frac{\lambda_u}{2} \int_{D_c} u^2 \,d\Omega + \frac{1}{2} \int_{D_a}{\bar{p}_{s} p_{s}	\,d\Omega}
	\\
	s.t. &
	\begin{dcases}
	\label{eq:state_dyn}
		-\nabla \cdot (a\nabla p_{s}) - b \omega^2 p_{s} = f		& \text{in } \Omega
		\\
		a \nabla p_{s} \cdot \mathbf{n} = g	  &	\text{on } \Gamma_i
		\\
		a\nabla p_{s}\cdot\mathbf{n} + \alpha p_{s}  =  0 		& \text{on } \Gamma_e
	\end{dcases}
\end{align}
where:
\begin{equation*}
    	\begin{dcases}
		f = \omega^2(b-b_0)p_{i} +\nabla\cdot[(a-a_0)\nabla p_{i}] \\
		g = - a \nabla p_i \cdot \mathbf{n}
		\\
		\alpha = a\Big(jk_0 + \frac{1}{2R}\Big)
	\end{dcases}
\end{equation*}
and the functional relationships between the control functions $u$ and $v$ and the perturbed material properties are:
\begin{equation*}
    a = a_0\,e^{-v} \qquad b = b_0\,e^{-u}
\end{equation*}
in this way the positivity of the density $\rho$ and bulk modulus $\kappa$ is ensured for any choice of the control functions $u$ and $v$. The exponential change of variables to ensure positivity of the control variables is standard and was used in \cite{chen2021optimal} when controlling the wave propagation velocity.

We now derive a set of first-order optimality conditions applying the Lagrangian method \cite{Trol2010}. Using this idea, we obtain an explicit expression for the gradient of the cost functional in the continuous setting. First of all, we define suitable functional spaces for state and control functions. We select the complex-valued Hilbert space $H^1(\Omega)$ as the state space, that is $\mathcal{V}=H^1(\Omega)$. The state problem is well-posed as long as its coefficients $\kappa$ and $\rho$ are bounded and positive \cite{colton1998inverse}. Since we have selected an exponential modulation of background properties we can select as control space $\mathcal{U}=L^{\infty}(D_c)^2$, that is the space of real-valued two-dimensional vector functions which are essentially bounded. In other words, for each $\bb{x} \in D_c$ we associate a real-valued control pair $(u(\bb{x} ),v(\bb{x} ))$ whose elements are bounded.

\noindent The Lagrangian functional $\mathcal{L}\colon \mathcal{V} \times \mathcal{U} \times \mathcal{W}^{*} \to \R$ can be formed as:
\begin{equation}\label{eq:Last L}
	\LL \coloneqq J
	+ \Re\Big\{\int_{\Omega} (\nabla \cdot (a\nabla p_{s}) + b \omega^2 p_{s} + f)\bar\lambda \, d\Omega\Big\}
\end{equation}
where the adjoint function $\lambda\colon\Omega \to \C$ belongs to $H^1(\Omega)$, that is we can identify $\mathcal{W}^{*} = H^{1}(\Omega)$. Note that the Lagrangian is defined as a real-valued functional and an equivalent formulation can be recovered by using the imaginary part.

\noindent A system of first-order necessary conditions for optimality is obtained by taking the G\^ateaux derivatives of the Lagrangian with respect to state, control and adjoint variables independently (see e.g.\ \cite{Trol2010}). The adjoint dynamics is obtained by setting to zero the Lagrangian derivative with respect to an arbitrary state variation $\phi \in H^1(\Omega)$. Applying the divergence theorem and substituting the boundary conditions, the Lagrangian can be rewritten as:
\begin{equation}
	\LL =  J
	+ \Re\Big\{\int_{\Omega} -a\nabla p_{s} \cdot \nabla \bar\lambda  + b \omega^2 p_{s} \bar\lambda + f\bar\lambda \,d\Omega + \int_{\Gamma_i} g \bar\lambda \,d\Gamma + \int_{\Gamma_e} \alpha p_s \bar\lambda \,d\Gamma\Big\}.
\end{equation}

\noindent 
$\mathcal{L}$ is a functional which maps complex-valued functions to real numbers, therefore to compute its G\^ateaux derivatives we make use of basic results from complex analysis, that is we apply Wirtinger's calculus rules \cite{wirtinger1927formalen}, in particular recall that  $\frac{d\Re\{ c\,z  \}}{dz}=c/2$ and $\frac{d(\bar{z}\,z)}{dz}=\bar z$ for $c\,,\,z\, \in \C$.

\noindent Hence, the G\^ateaux derivative of $\mathcal{L}$ with respect to $p_{s}$ is:
\begin{equation}
\begin{aligned}
	\LL_{p_{s}}'[\varphi] = 
	\frac{1}{2}\int_{D_a}{\bar{p}_{s} \varphi \,d\Omega}
	+ \frac{1}{2}\int_{\Omega}{- a \nabla \bar\lambda \cdot \nabla \phi  + b \omega^2  \bar\lambda \phi \, d\Omega} + \frac{1}{2}\int_{\Gamma_e} \alpha \, \bar\lambda \, \phi \, d\Gamma = 0 \quad \forall \phi \in H^1(\Omega)
\end{aligned}
\end{equation}
which is the weak formulation of the adjoint dynamics:
\begin{equation}\label{eq:PDE adjoint}
	\begin{array}{ll}
		-\nabla \cdot (a\nabla \lambda) - b\omega^2 \lambda = p_{s}\, \chi_{D_a}		& \text{in } \Omega \\
		\nabla \lambda \cdot \mathbf{n} = 0	  &	\text{on } \Gamma_i \\
		\nabla \lambda\cdot\mathbf{n} + \Big( \frac{1}{2R} - j k_0 \Big) \lambda  =  0 		& \text{on } \Gamma_e
	\end{array}
\end{equation}
being $\chi_{D_a}(\bb{x})$ the indicator function of the domain $D_a$. 

We now turn to the optimality conditions involving the control functions $u$ and $v$. The Lagrangian \eqref{eq:Last L} can be rewritten substituting the explicit form of $f$ and considering satisfied the boundary conditions of the state PDE as:

\begin{equation}\label{eq:control_Lag}
	\LL = J
	+ \Re\Big\{\int_\Omega{b\omega^2(p_{s} + p_{i})\bar\lambda	\,d\Omega} 
	-\int_\Omega{a\nabla (p_{s}+p_{i})\cdot\nabla \bar\lambda	\,d\Omega}
	+ \int_{\Gamma_e}{a\nabla (p_{s}+p_{i})\cdot\mathbf{n}\,\bar\lambda	\,d\Gamma}\Big\}
\end{equation}
so that we can easily take control variations $\psi \in L^{\infty}(D_c)$  Physically, the control variations cannot modify the background properties outside of the cloak. Hence,  $a=a_0$ on $\Gamma_e$, $a_v'[\psi] = b_u'[\psi] = 0$ on $D_a$ and
\begin{equation*}
	a_v'[\psi]=\{a_0e^{-v}\}_v'[\psi] = -a_0e^{-v}\psi = -a\psi
\end{equation*}
and similarly $b_u'[\psi] = -b\psi$, so that the control necessary conditions (i.e.\ the reduced gradient) in variational form results in: 
\begin{equation}
	\begin{dcases}
		\LL_v'[\psi] = \lambda_v\int_{D_c}{v\psi	\,d\Omega} 
		+ \Re\Big\{\int_{D_c}{a\nabla(p_{s} + p_{i})\cdot\nabla \bar\lambda\,\psi	\,d\Omega}\Big\} = 0 &\quad \forall \psi \in L^{\infty}(D_c)
		\\
		\LL_u'[\psi] = \lambda_u\int_{D_c}{u\psi	\,d\Omega} 
		- \Re\Big\{\int_{D_c}{b\omega^2(p_{s} + p_{i}) \bar\lambda \,\psi	\,d\Omega}\Big\}  = 0 &\quad \forall \psi \in L^{\infty}(D_c).
	\end{dcases}
\end{equation}
The strong form of the reduced gradient can be identified as:

\begin{equation}
\label{red_grad}
\begin{aligned}
&\nabla J_v = \lambda_v v + \Re\big\{a\nabla(p_{s} + p_{i})\cdot\nabla \bar\lambda\big\} \\
&\nabla J_u = \lambda_u u -  \Re\big\{b\omega^2 (p_{s} + p_{i}) \bar\lambda\big\}. \\
\end{aligned}
\end{equation}

\begin{figure}
	\centering
	\includegraphics[width=0.45\textwidth]{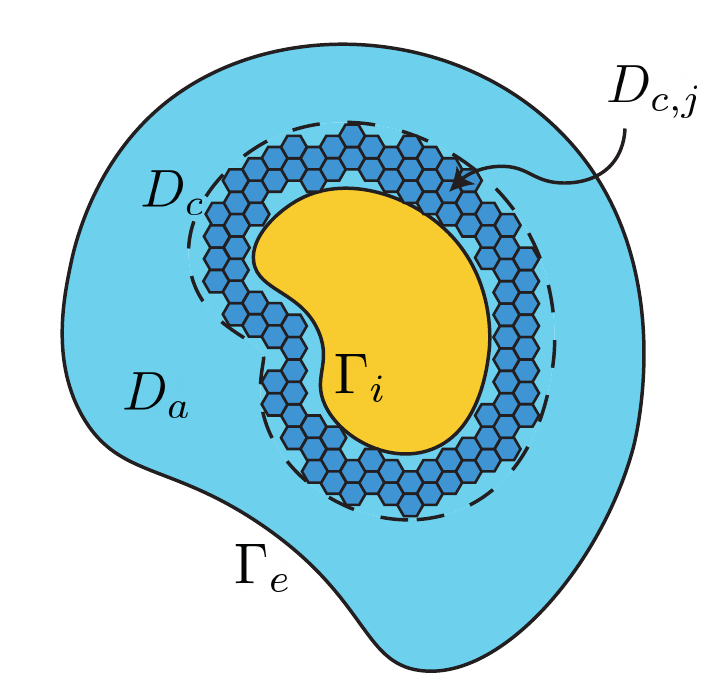}
	\caption{the obstacle is surrounded by hexagonal domains $D_{c,j}$. Notice that $\cup_{j=1}^{N_c} D_{c,j}\,\subset\,D_c$ is the control domain.}
	\label{fig:psi definition}
\end{figure}

\noindent Equations \eqref{red_grad} together with the adjoint Equation \eqref{eq:PDE adjoint} and the state Equation \eqref{eq:state} constitute a system of first-order necessary conditions for optimality. Note that we did not make any assumption on the structure of the control basis functions other than belonging to the space $L^{\infty}(D_c)$. However, the actual controlled material properties will be realized with piece-wise constant functions at the microstructure level. In order to preserve the optimal properties at the microstructure, we express the control functions as linear combinations of indicator functions describing the cell domain. In particular, let us define a subdivision of the control domain $D_c$ in $N_c$ disjoint sets whose elements $D_{c,j}$ satisfy:
\begin{equation*}
\bigcup_{j=1}^{N_c} D_{c,j} \, \subseteq \, D_c \quad \textrm{and} \quad D_{c,j} \cap D_{c,i} = \emptyset \quad \textrm{for} \quad  i \neq j
\end{equation*}
and define the functions $\psi_{j}(\mathbf{x}) = \chi_{D_{c,j}}(\mathbf{x})$ as the indicator functions of such sets. Then it is natural to express to control variables $u$ and $v$ as:
\begin{equation}
\label{semi_disc}
    u = \sum_{j=1}^{N_c} \psi_j(\mathbf{x}) u_j =  \boldsymbol{\psi}(\mathbf{x})^{\top} \mathbf{u} \qquad v = \sum_{j=1}^{N_c} \psi_j(\mathbf{x}) v_j = \boldsymbol{\psi}(\mathbf{x})^{\top} \mathbf{v}
\end{equation}
where the shape functions $\psi_j$ are defined according to the cell shape and distribution in the domain $D_c$ and the constant coefficients $u_i$ and $v_i$ of the linear combination are the control variables of the optimization problem.  The control discretization layout is shown in Figure~\ref{fig:psi definition}. This formulation allows to preserve cloak's optimal properties at the microstructure level, as was mentioned in the introduction.

The optimization problem is still set in the infinite-dimensional space for the state and adjoint variables. It is also clear that $u = \sum_{j=1}^{N_c} \psi_j(\mathbf{x}) u_j \in L^{\infty}(D_c)$ and the functional setting of the OCP is still consistent. Regarding the optimality conditions, it is easy to see that state and adjoint dynamics are unchanged. Slightly more care is needed to recover the form of the reduced gradients $\nabla J_{v_j}$ and $\nabla J_{u_j}$ for $j=1,\ldots,N_c$. We substitute Equations~\eqref{semi_disc} in the Lagrangian formulation \eqref{eq:control_Lag} as:
\begin{equation*}
\begin{aligned}
& \mathcal{L} = \frac{\lambda_v}{2} \mathbf{v}^{\top} \Big(\int_{D_c} \boldsymbol{\psi} \boldsymbol{\psi}^{\top} \, d\Omega\Big) \mathbf{v} + \frac{\lambda_u}{2} \mathbf{u}^{\top} \Big(\int_{D_c} \boldsymbol{\psi} \boldsymbol{\psi}^{\top} \, d\Omega\Big) \mathbf{u} + \frac{1}{2}\Re\Big\{ \int_{D_a}{\bar{p}_{s} p_{s}	\,d\Omega}\Big\}\\ 
&+\Re\Big\{\int_\Omega{\omega^2b(\mathbf{u})(p_{s} + p_{i})\bar\lambda	\,d\Omega}
	-\int_\Omega{a(\mathbf{v})\nabla (p_{s}+p_{i})\cdot\nabla \bar\lambda	\,d\Omega}
	+ \int_{\Gamma_e}{a(\mathbf{v})\nabla (p_{s}+p_{i})\cdot\mathbf{n}\,\bar\lambda	\,d\Gamma}\Big\}
\end{aligned}
\end{equation*}
where :
\begin{equation*}
    a(\bb{v}) = a_0\,e^{-\bs{\psi}^{\top}\bb{v}} \qquad b(\bb{u}) = b_0\,e^{-\bs{\psi}^{\top}\bb{v}}.
\end{equation*}
Furthermore, since $\forall \bb{x} \in D_c$ there is at most one index $k$ such that $\psi_k(\bb{x})\neq 0$ we  have:
\begin{equation}
\label{nice_psi}
    e^{-\bs{\psi}^{\top}\bb{v}} =  e^{-\sum_{j=1}^{N_c} v_j\,\psi_j}  = \sum_{j=1}^{N_c} e^{-v_j} \psi_j
\end{equation}
for every vector $\mathbf{v} \in \R^{N_c}$. Note also that the gradient of $a$ and $b$ can be written as:
\begin{equation*}
    \nabla_{\bb{v}} a = - \bs{\psi} \, a \qquad \nabla_{\bb{u}} b = - \bs{\psi} \, b 
\end{equation*}
so that the reduced gradients can be expressed by taking the finite-dimensional gradient of the Lagrangian with respect to $\bb{v}$ and $\bb{u}$, that is:
\begin{equation}
\label{red_grad_disc_1}
\begin{aligned}
\nabla J_{\mathbf{v}} &= \lambda_v \Big(\int_{D_c} \bs{\psi} \bs{\psi}^{\top} \, d\Omega\Big) \mathbf{v} +  \Re\Big\{\int_{D_c} \bs{\psi} \, a \nabla(p_{s} + p_{i})\cdot\nabla \bar\lambda \, d\Omega\Big\}
\\
\nabla J_{\mathbf{u}} &= \lambda_u \Big(\int_{D_c} \boldsymbol{\psi} \boldsymbol{\psi}^{\top} \, d\Omega\Big) \mathbf{u} + \Re\Big\{\int_{D_c} \boldsymbol{\psi} \, b\omega^2 (p_{s} + p_{i}) \bar\lambda \, d\Omega \Big\}.\\
\end{aligned}
\end{equation}
Note that $\displaystyle \int_{D_c} \boldsymbol{\psi} \boldsymbol{\psi}^{\top} \, d\Omega$ is a diagonal matrix whose entries are the areas of the respective cells. We can now turn to the full discretization of the problem.

\subsection*{Discretization of the OCP}
For  the  numerical  solution  of the  OCP  we  employ  the  Finite  Element Method (FEM). We select piecewise quadratic, globally continuous ansatz functions $\phi_i$ ($\mathbb{P}_2$ finite elements) for the space approximation of state and adjoint in $\Omega$ while the control basis functions do not need any spatial approximation since their functional form is expressed by Equation \eqref{semi_disc}. The FEM approximation of the state equation reads:
\begin{equation*}
    A(\mathbf{u},\mathbf{v}) \mathbf{p} = \mathbf{f}(\mathbf{u},\mathbf{v})
\end{equation*}
where:
\begin{equation*}
    A_{ij} = \int_{\Omega} a(\mathbf{v}) \nabla \phi_i \cdot \nabla \phi_j - \int_{\Omega} \omega^2 b(\mathbf{u}) \phi_i \phi_j -\int_{\Gamma_e} \alpha \phi_i \phi_j d\Gamma 
\end{equation*}
since $a=a_0$ in $\Omega \setminus D_c$ it is useful to rewrite:
\begin{equation*}
    a = a_0 + a_0( e^{-\boldsymbol{\psi}^{\top}\mathbf{v}}-1) \qquad b = b_0 + b_0( e^{-\boldsymbol{\psi}^{\top}\mathbf{u}}-1)
\end{equation*}
so that using Equation \eqref{nice_psi} the components of A can be separated as:
\begin{equation}
\label{A_structure}
\begin{aligned}
    & A_{ij} = \int_{\Omega} a_0 \nabla \phi_i \cdot \nabla \phi_j d\Omega + \sum_{k=1}^{N_c} (e^{-v_k}-1) \int_{D_{c,k}} a_0\nabla \phi_i \cdot \nabla \phi_j d\Omega - \int_{\Omega} b_0  \omega^2 \phi_i \phi_j d\Omega \\
    & - \sum_{k=1}^{N_c} (e^{-u_k}-1) \int_{D_{c,k}} b_0 \omega^2 \phi_i \phi_j d\Omega -\int_{\Gamma_e} a_0(jk_0 + \frac{1}{2R}) \phi_i \phi_j d\Gamma \\
    & = \big(A_0\big)_{ij} + \sum_{k=1}^{N_c} (e^{-v_k}-1) \big(A_k\big)_{ij} + \big(B_0\big)_{ij}  + \sum_{k=1}^{N_c} (e^{-u_k}-1) \big(B_k\big)_{ij} + C
\end{aligned}
\end{equation}
where the matrices $A_{k}$,$B_{k}$ and $C$ can be precomputed and only their sum must be performed when varying the control vectors $\mathbf{u}$ and $\mathbf{v}$. Besides the presence of the exponential function that enforces the positive definiteness of the material properties, Equation \eqref{A_structure} highlights the bilinear structure of the control problem. Finally, it is easy to notice that the matrix A is symmetric being the sum of symmetric matrices. The components of the right-hand side $\mathbf{f}$ can be written as:
\begin{equation*}
\begin{aligned}
    f_i & = \int_{\Omega} (a-a_0) \nabla p_{i} \cdot \nabla \phi_i d\Omega + \int_{\Omega} (b-b_0) \omega^2 p_i \phi_i d\Omega - \int_{\Gamma_i} a_0 \nabla p_i \cdot \mathbf{n} \, \phi_i \, d\Omega \\
    & = \sum_{k=1}^{N_c} (e^{-v_k}-1) \int_{D_{c,k}} \nabla p_i \cdot \nabla \phi_i d\Omega + \sum_{k=1}^{N_c} (e^{-u_k}-1) \int_{D_{c,k}} b_0 \omega^2 p_i \phi_i d\Omega - \int_{\Gamma_i} a_0 \nabla p_i \cdot \mathbf{n}\,\phi_i \, d\Omega \\
    & = \sum_{k=1}^{N_c} (e^{-v_k}-1) \big(l_{k}\big)_i + \sum_{k=1}^{N_c} (e^{-u_k}-1) \big(d_{k}\big)_i - q_i
\end{aligned}
\end{equation*}
\begin{figure}
	\centering
	\includegraphics[width=0.8\textwidth]{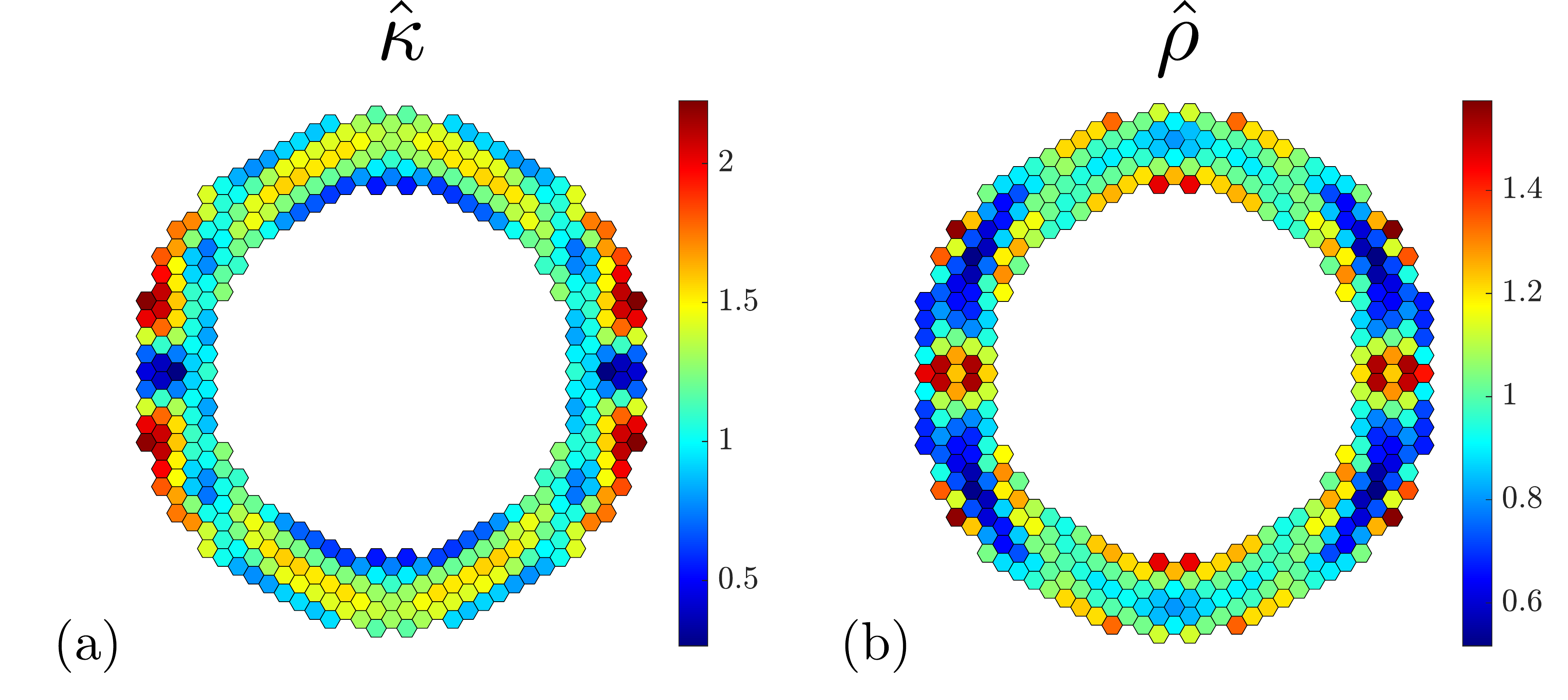}
	\caption{Optimal distribution of material properties computed solving the OCP when the shape of the target $\Gamma_i$ is a circle and the angular frequency $\omega$ is selected such that $\lambda/r=0.69$, $r$ being the radius of the target. The incidence direction is $\mathbf{a}=[1, 0]$. 390 hexagonal unit cells of edge $l =\SI{8.7}{\percent}\,\lambda$ are employed. The external radius of the cloak is $1.57\,r$.}
	\label{fig:unc_ocp}
\end{figure}
\begin{figure}
	\centering
	\includegraphics[width=0.8\textwidth]{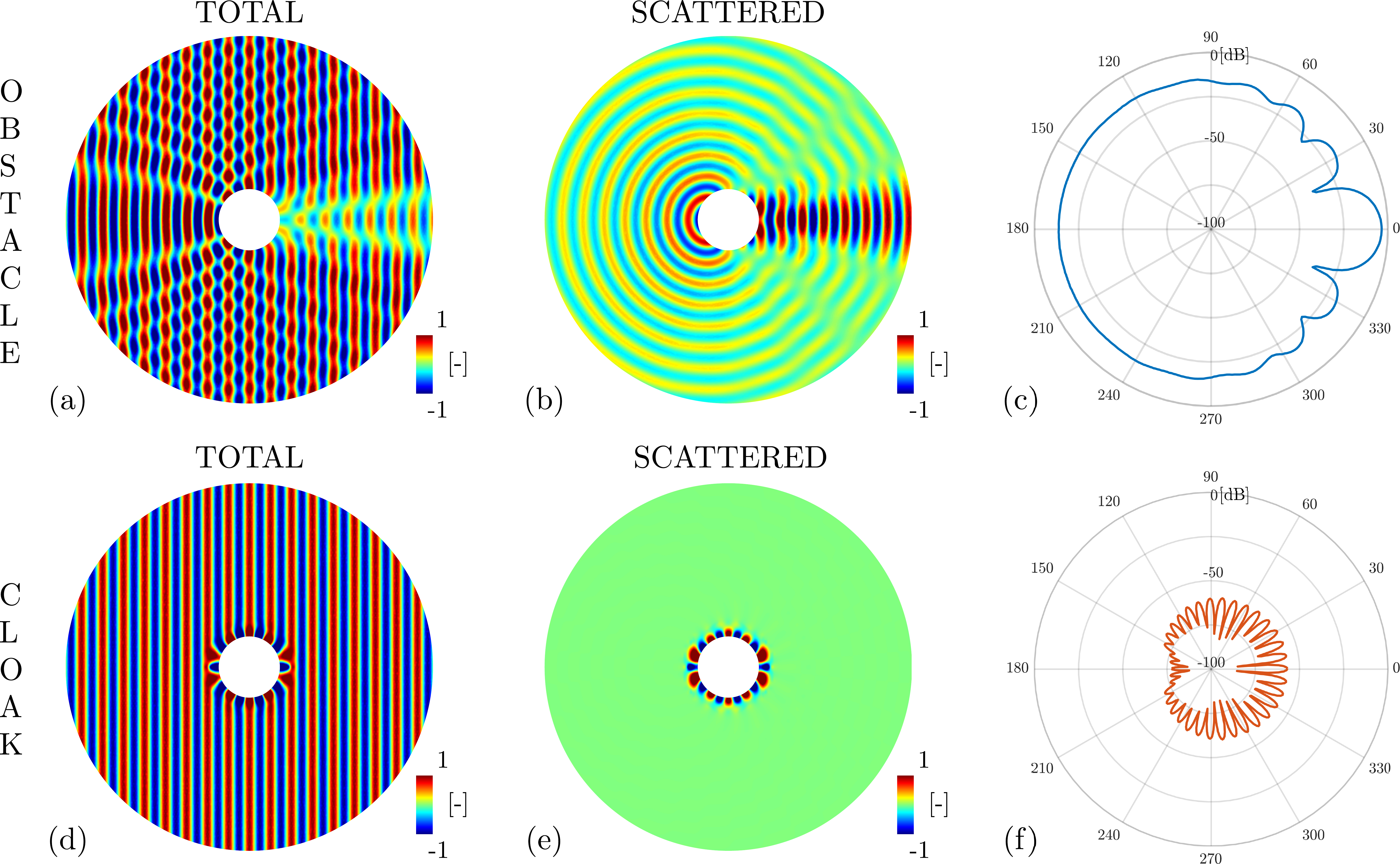}
	\caption{(a) Total pressure field, obstacle case. The pressure is normalized with respect to the amplitude of the incident wave. (b) Scattered pressure field, obstacle case. (c) Decibel reduction in acoustic intensity computed for the scattered field at 1 meter from the obstacle, with respect to the incident intensity. (d) Total pressure field, cloak case, when the material properties in the cloak are the optimal ones shown in Figure \ref{fig:unc_ocp}. (e) Scattered pressure field, cloak case. (f) Decibel reduction in acoustic intensity with respect to the incident one, cloak case. }
	\label{fig:scattering}
\end{figure}
which again shows the same bilinear structure in the way the control functions enter the right-hand side. The adjoint discretization follows the same steps for the left-hand side while the right-hand side corresponds to the FEM discretization of the of the state projected in $L^2(\Omega)$. That is we have:
\begin{equation*}
    A^\dag(\bb{u},\bb{v}) \bs{\lambda} = M_{D_a} \bb{p}
\end{equation*}
where $M_{D_a}$ is the restriction of the usual mass matrix to the observation domain $D_a$, that is the domain in which we want to minimize the scattered field, and $(\cdot)^\dag$ is the Hermitian operator. Note that since $A^\dag = \bar A$ the discretized version of the state operator remains self-adjoint. Finally, the FEM discretization of state $\mathbf{p}$ and adjoint $\boldsymbol{\lambda}$ can be plugged in Equation \eqref{red_grad_disc_1} to obtain the fully discrete version of the reduced gradient, that at component level of $v_k$ can be written as:
\begin{equation}
\label{red_grad_disc_2}
\begin{aligned}
\nabla J_{v_k} &= \lambda_v |D_{c,k}|v_k +  \,e^{-v_k}  \int_{D_{c,k}} \, a_0 \, \nabla(p_{s} + p_{i})\cdot\nabla \bar\lambda \, d\Omega  \\
& = \lambda_v |D_{c,k}|v_k +  \,e^{-v_k}  \bs{\lambda}^\dag\Big( A_{k} \bb{p} + \bb{l}_{k} \Big) \\
\end{aligned}
\end{equation}
where $|D_{c,k}|$ is the measure of the set associated to the $k$\textsuperscript{th} cell. For the $u_k$ control vectors we have:
\begin{equation*}
\begin{aligned}
\nabla J_{u_k} &= \lambda_u |D_{c,k}|u_k + e^{-u_k} \int_{D_{c,k}}  b_0\omega^2 (p_{s} + p_{i}) \bar\lambda \, d\Omega \\
&= \lambda_u |D_{c,k}|u_k +e^{-u_k}\bs{\lambda}^\dag\Big( B_{k}\bb{p} + \bb{d}_{k} \Big).
\end{aligned}
\end{equation*}
Once the fully discretized version of the optimality conditions is obtained, we setup Algorithm~\ref{alg_steep} using an iterative steepest descent method to solve the OCP with microstructure specified by the functions $\psi_k$.

\begin{algorithm}
  \begin{algorithmic}[1]
  \State
  \State $\boldsymbol{\psi},N_c \gets \textrm{Define microstructure shape and domain}$
  \State $\textrm{FEM Model} \gets \text{Assemble constant FEM matrices} $
  \State $\mathbf{u}^{0},\mathbf{v}^{0} \gets \text{Assign control initial
  guesses}$ 
  \State
    \For{$t=1:\text{maxIter}$} 
        \State
        \State    $\bb{p}^t \gets \text{Solve state equation: \,\,\,\,} A(\mathbf{u}^{t},\mathbf{v}^{t})\mathbf{p} = \mathbf{f}(\mathbf{u}^t,\mathbf{v}^t)$ 
        \State
        \State    $\bs{\lambda}^t \gets \text{Solve adjoint equation: \,\,} \bar{A}(\mathbf{u}^{t},\mathbf{v}^{t})\bs{\lambda} = M_{D_a}\bb{p}^t $ 
        \State
        \State $\nabla J(\bb{u}^{t})_k \gets \lambda_u |D_{c,k}|u_k^t +e^{-u^t_k}\bs{\lambda}^{t\,\dag}\Big(B_{k}\bb{p}^{t} + \bb{d}_{k} \Big)$ 
        \State 
        \State $\nabla J(\bb{v}^{t})_k \gets  \lambda_v |D_{c,k}|v^t_k +  \,e^{-v^t_k}  \bs{\lambda}^{t\,\dag}\Big( A_{k}\bb{p}^{t} + \bb{l}_{k} \Big)$ 
        \State
        \If{ $\norm{\Big( \nabla J(\bb{u}^{t})\, , \, \nabla J(\bb{v}^{t}) \Big)} < \text{tol}$}
        \State $  \text{return}$
        \EndIf
        \State
        \State $\tau \gets \text{ArmijoBacktracking}(J,\nabla J(\bb{u}^{t}),\nabla J(\bb{v}^{t}),\bb{u}^t,\bb{v}^t) $ 
        \State
        \State $\bb{u}^{t+1} \gets \bb{u}^{t} - \tau\nabla J(\bb{u}^{t})  $
        \State $\bb{v}^{t+1} \gets \bb{v}^{t} - \tau\nabla J(\bb{v}^{t})  $  \Comment{Update control}

        \State
    \EndFor
    
  \end{algorithmic}
  \caption{Steepest Descent for Optimal Cloak }
  \label{alg_steep}
\end{algorithm}
The OCP is solved for a circular target surrounded by the set of $N_c=390$ hexagonal unit cells as shown in Figure~\ref{fig:unc_ocp} when probed by acoustic illumination from left to right at an angular frequency $\omega$ corresponding to $\lambda/r=0.69$, with $\lambda=2\pi c_0/\omega$ and $r$ being the radius of the target. The solution is shown in Figure~\ref{fig:unc_ocp} in terms of nondimensional material properties $\hat{\rho}=\rho/\rho_0$ and $\hat{\kappa}=\kappa/\kappa_0$. Figure~\ref{fig:scattering} compares the total (Figures~\ref{fig:scattering}(a)~vs~(d)) and scattered pressure fields (Figures~\ref{fig:scattering}(b)~vs~(e))  between the cloaked and uncloaked case. The mean scattered intensity $I_{mean}$ at one meter from the surface of the obstacle is also computed for all the azimuthal angles $\theta$ and in Figures~\ref{fig:scattering}(c) and (f) is shown for comparison in terms of Decibel reduction with respect to the incident intensity $I_{inc}$:
\begin{equation*}
    \Delta(\theta)=10 \log_{10}\left( \frac{I_{mean}(\theta)}{I_{inc}} \right) \quad [\operatorname{dB}].
\end{equation*}
An average reduction of \SI{65}{[\decibel]} on the scattered intensity is obtained using the hexagonal microstructure discretization.  Note that, as shown in Figure~\ref{fig:unc_ocp}, we need to obtain equivalent controlled properties both higher and smaller than those of the background fluid.

\section{Unit Cell Design}
\begin{figure}
	\centering
	\includegraphics[width=0.8\textwidth]{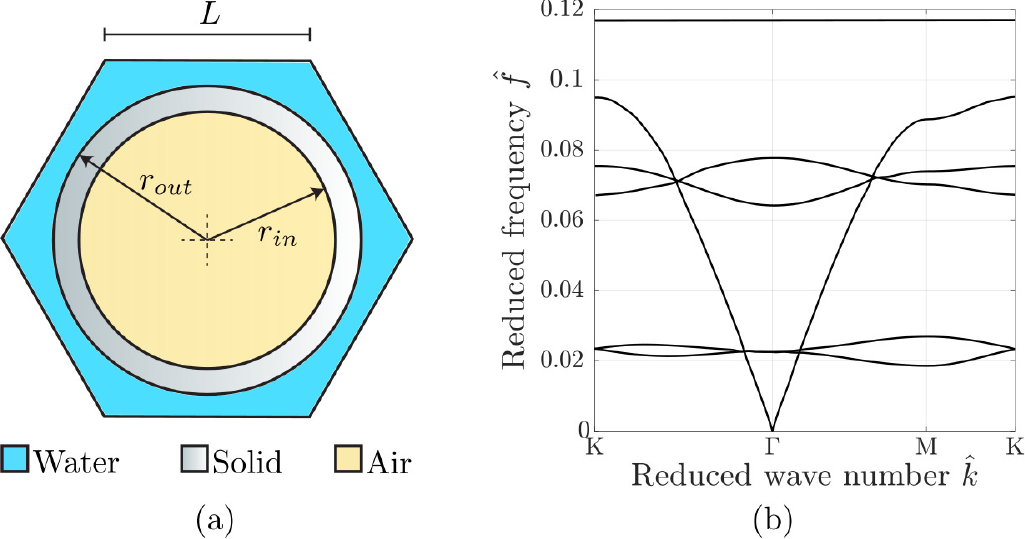}
	\caption[]{(a) Schematic of a unit cell comprising a circular inclusion made by a solid phase filled with air. (b)~Dispersion relation computed along the boundary of the Irreducible Brillouin Zone. The reduced frequency is computed as $\hat{f}=f L/c_0$, while the reduced wavenumber stands for the adimensional $kL$, which spans between the high symmetry points $\Upgamma$, K and M.}
	\label{fig:cell1}
\end{figure}

The required material parameter distribution obtained through the solution of the OCP introduced in the previous section has to be practically realized with opportunely designed microstructures that show the appropriate equivalent density and bulk modulus when homogenized. It is well known \cite{laude2015phononic}, that hexagonal lattices of solid inclusions in water behave in the long-wavelength limit as isotropic acoustic fluids, whose properties can be tailored upon control on the material and shape of the inclusion itself. For this very reason the cloak sub-domains have been chosen to be shaped as hexagons: in this way they can naturally be filled by hexagonal lattices. The basic configuration considered in the bidimensional setting consists thus of a circular inclusion placed in each lattice point and made by a material with high contrast with respect to the hosting medium, e.g.\ a metal. This allows to obtain a wide range of material properties with densities and bulk moduli that are generally higher than that of water. Preliminary results shown in the previous section (Figure \ref{fig:unc_ocp}) underline the need to go also for $\rho$ and $\kappa$ smaller than those of water: it is thus implied that some kind of porosity has to be contemplated in the solid inclusion. Indeed, since resonance phenomena are not exploited in this application, the density can be simply evaluated with the rule of mixtures:
\begin{equation}
    \rho_{\hom}=\sum_i \chi_i\rho_i
\end{equation}
where $\rho_i$ is the density of the constituents and $\chi_i$ is the cell volume filling fraction of each constituent. This in turn implies that a third light phase has to be included in the mix other than the fluid and the solid. The simplest configuration considered consists thus of a hollow cylinder filled by air ($\rho=\SI{1.23}{ [\kilo\gram\per\cubic\meter]}$, $\kappa=\SI{0.14}{[\mega\pascal]}$, ref Figure~\ref{fig:cell1}(a)).\\ 
The equivalent bulk modulus is instead computed via inspection of the dispersion relation of each considered lattice, computed via Bloch analysis on the unit cell \cite{laude2015phononic}. A typical dispersion relation is shown in Figure~\ref{fig:cell1}(b): in the long wavelength limit, the linearity of the branch justifies the evaluation of $\kappa_{\hom}$ as:
\begin{equation}
    \kappa_{\hom}=c_{ph}^2\rho_{\hom}
\end{equation}
where $c_{ph}$ is the phase speed computed as the slope of the very branch emanating from the origin. In order to compute the set of obtainable $\rho_{\hom}$, $\kappa_{\hom}$, the geometry is parametrized with the two characteristic adimensional parameters $\hat{r}_{out}=r_{out}/L$ and $\hat{r}_{in}=r_{in}/L$ (Figure~\ref{fig:cell1}(a)), whose variation is considered to be bounded in the following way:
\begin{equation}
    \begin{cases}
    \hat{r}_{in}\geq\delta \hat{r}_1\\
    \hat{r}_{out}\leq\frac{\sqrt{3}}{2}-\delta \hat{r}_1\\
    \hat{r}_{out}\geq\hat{r}_{in}+\delta \hat{r}_2
    \end{cases}
\end{equation}

\begin{figure}
	\centering
	\includegraphics[width=0.8\textwidth]{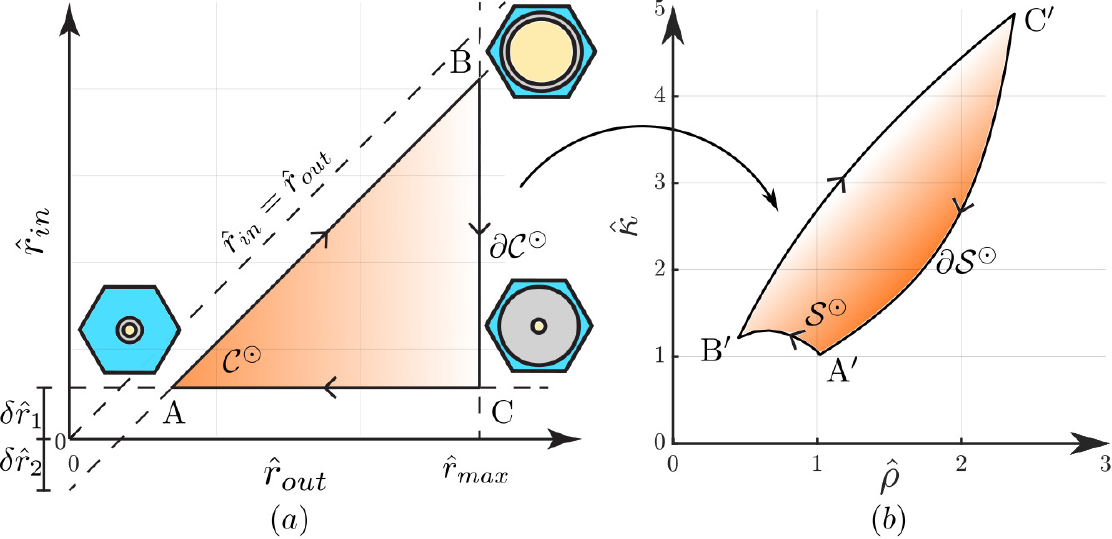}
	\caption[]{(a) The feasible set $\mathcal{C}^\odot$ of unit cell geometries comprising a cylindrical inclusion defined in the $\hat{r}_{out} \times \hat{r}_{in}$ space, which is enclosed by the curve $\partial\mathcal{C^\odot}$. (b) The corresponding reachable set of material properties $\mathcal{S}^\odot$ is enclosed in the curve $\partial\mathcal{S^\odot}$, obtained performing the long-wavelength homogenization of the parametrized cell along curve $\partial\mathcal{C^\odot}$. }
	\label{fig:cell2}
\end{figure}
\begin{figure}[t]
	\centering
	\includegraphics[width=0.3\textwidth]{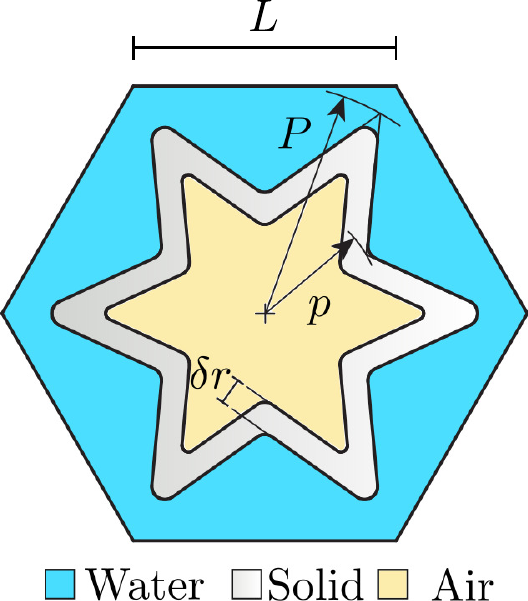}
	\caption[Domain definition]{Schematic of the unit cell comprising a hollow N-pointed star.}
	\label{fig:cell3}
\end{figure}
\begin{figure}[t]
	\centering
	\includegraphics[width=0.9\textwidth]{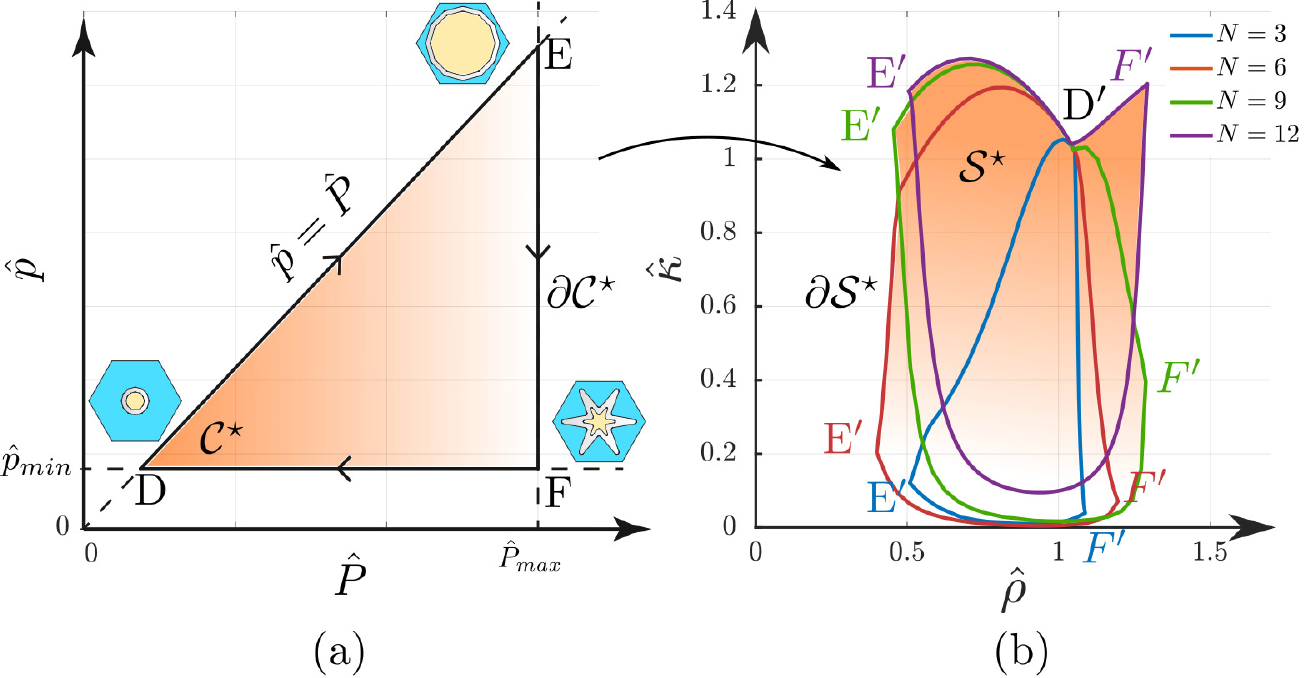}
	\caption[]{(a) The feasible set of unit cell geometries $\mathcal{C^\star}$ comprising a N-pointed star inclusion defined in the $\hat{P} \times \hat{p}$ space, which is enclosed by the curve $\partial\mathcal{C^\star}$. (b) The corresponding reachable set of material properties $\mathcal{S^\star}$ is enclosed in the curve $\partial\mathcal{S^\star}$, obtained performing the long-wavelength homogenization for each cell described by the points on the curve $\partial\mathcal{C^\star}$.}
	\label{fig:cell4}
\end{figure}


\noindent where $\delta \hat{r}_{1,2}$ are the nondimensional minimum  feature sizes, that is the thinnest gap and wall allowed. These constraints define a closed feasible region $\mathcal{C}^\odot$ in the plane $\hat{r}_{in} \times \hat{r}_{out}$, that is shown in Figure~\ref{fig:cell2}(a). By computing the homogenized properties of the associated lattices, the contour $\partial\mathcal{C^\odot}$ going across the extremal points ABC is mapped to a curve $\partial\mathcal{S^\odot}$ joining A$^\prime$B$^\prime$C$^\prime$ in the $\hat{\rho} \times \hat{\kappa}$ space, with $\hat{\rho}=\rho_{\hom}/\rho_0$ and $\hat{\kappa}=\kappa_{\hom}/\kappa_0$. This defines the set of the obtainable material properties $\mathcal{S^\odot}$. In Figure~\ref{fig:cell2}(b) such curve is computed for a configuration where the solid phase is chosen to be aluminium ($\rho=\SI{2700}{[\kilo\gram/\cubic\meter]}$, Young's Modulus $YM=\SI{70}{[\giga\pascal]}$, Poisson's ratio $\nu=0.3$) and the minimum features are selected as $\delta \hat{r}_2=\SI{5}{\percent}$ and $\delta \hat{r}_1=4\delta \hat{r}_2$. It can be seen how the inclusion of the light phase allows for obtaining $\hat{\rho}<1$, notice however how it is hard to reach the region where $\hat{\kappa}$ is less than $1$. To enlarge the feasible region, another configuration is thus considered: the inclusion is now shaped as a N-pointed star, N being a multiple of $3$; other than maintaining the invariance of the lattice upon rotation of $\pi/3$, i.e.\ the symmetry required for isotropy, the oblique walls allow to reduce the tangential stiffness of the inclusion. When considering hydrostatic loads, this in turn increases the compressibility with respect to the case of the hollow cylinder. A N-pointed star is completely characterized by the lengths of the internal and external tips $\hat{P}=P/L$ and $\hat{p}=p/L$, by the fillet radii and by the thickness of the wall (Figure~\ref{fig:cell3}(a)). The latter two parameters are considered fixed and are chosen to be \SI{2.5}{\percent} and \SI{5}{\percent}, respectively. The bounds on the remaining two geometrical features are:
\begin{equation}
    \begin{cases}
    \hat{p}\geq\hat{p}_{min}\\
    \hat{p} \leq \hat{P}\\
    \hat{P} \leq \hat{P}_{max}
    \end{cases}
\end{equation}
these also define a feasible $\hat{p} \times \hat{P}$ region $\mathcal{C}^\star$ (Figure~\ref{fig:cell4}(a)) whose boundary $\partial\mathcal{C}^\star$ can be mapped to a path $\partial\mathcal{S}^\star$ in the $\hat{\rho} \times \hat{\kappa}$ space. In Figure~\ref{fig:cell4}(b) it is shown how adopting this type of unit cell the feasible set of material properties is enlarged also in the region that is not reachable with the hollow cylinder.
\begin{figure}[t]
	\centering
	\includegraphics[width=0.7\textwidth]{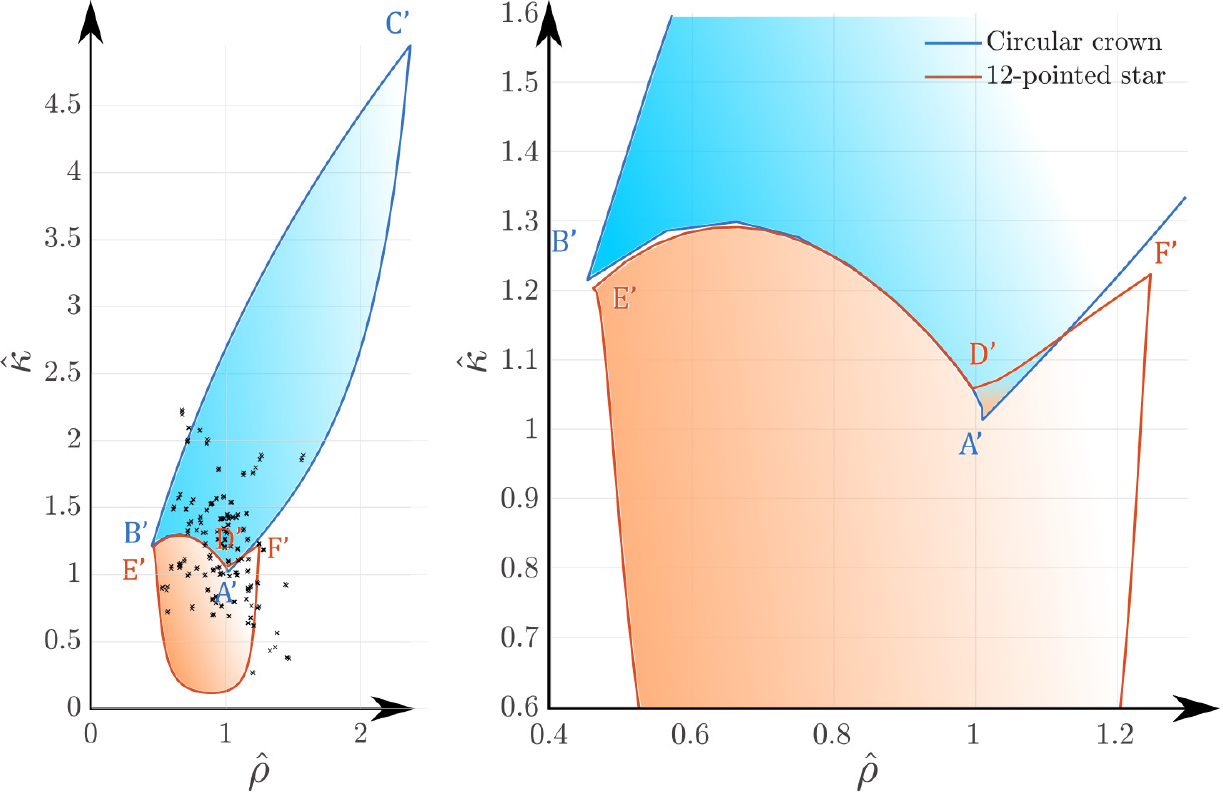}
	\caption[]{The overall reachable set of material properties given as the union of those obtained separately with the hollow cylinder inclusion and the 12-pointed star inclusion. A magnification around the material properties of the background fluid shows how the A$^\prime$B$^\prime$ and D$^\prime$E$^\prime$ curves almost overlap creating a connected set. Black markers in the graph are used to underline the location of the optimal material properties computed with the unconstrained OCP, that are shown in Figure~\ref{fig:unc_ocp}.}
	\label{fig:cell5}
\end{figure}

Note that, the higher the number N, the more similar is the N-pointed star to a hollow cylinder when $\hat{p} \rightarrow \hat{P}$. For this reason, the D$^\prime$E$^\prime$ curve for a 12-pointed star almost overlap with the A$^\prime$B$^\prime$ curve of the circular inclusion. This allows to obtain a connected feasible set $\mathcal{S}\coloneqq\mathcal{S}^\odot\cup\mathcal{S}^\star$ in the $\hat{\rho} \times \hat{\kappa}$ space, as shown in Figure~\ref{fig:cell5}, that will be considered in the following the reachable region for the equivalent material properties.

\section{Constrained Optimal Control Problem}
In this section we reformulate the fully discrete PDE-constrained optimization problem in order to satisfy the constraints imposed by the realization of the actual microstructure. That is we solve a reduced constrained optimization problem where the constrained control region generates equivalent material properties that lie in the reachable region $\mathcal{S}$ of the $\hat{\rho} \times \hat{\kappa}$ space. Furthermore, we include a regularization term in the control weightings to impose a smoother transition of material properties between neighboring cells.
First of all, the optimal material properties obtained in Section~\ref{sec:OCP} are plotted in Figure~\ref{fig:cell5} as black markers in the $\hat{\rho}  \times \hat{\kappa}$ plane. It can be noticed how part of them falls outside of the set of material properties that can be practically implemented by means of the microstructures described in the previous section. 
 
In order to constrain the control variables to lie on the feasible set described by the region $\mathcal{S}$, we equip the steepest descent Algorithm~\ref{alg_steep} with an additional projection step thus employing a standard Projected Gradient (PG) method \cite{nocedal2006numerical}. For each component-wise control pair $(v_k,u_k)$, the corresponding point $P_k=(\hat\rho_k,\hat\kappa_k)=(e^{v_k},e^{u_k})$ must lie in the region of the $\hat\rho\times\hat\kappa$ plane defined by $\mathcal{S}$. 

\noindent The feasible region in the control space is defined as 
$\mathcal{S}' = \{ (v,u) \in \R^{2} : (e^{v},e^{u}) \in \mathcal{S} \subset \R^{2} \}$  and we denote the projection onto $\mathcal{S}'$ as  $\Pi_{\mathcal{S}'}$. The pairwise vector projection  $\boldsymbol{\Pi}_{\mathcal{S}'}$ is defined as:
\begin{equation*}
 \Big(\boldsymbol{\Pi}_{\mathcal{S}'}(\bb{v},\bb{u})\Big)_k = \Pi_{\mathcal{S}'}(v_k,u_k)
\end{equation*}

\noindent The PG method consists of replacing the gradient update in Algorithm~\ref{alg_steep} with:
\begin{equation*}
(\bb{v},\bb{u})^{t+1}
=
\boldsymbol{\Pi}_{\mathcal{S}'}(\bb{v}^{t}-\tau \nabla J(\bb{v}^{t}) ,\bb{u}^{t}-\tau \nabla J(\bb{u}^{t}))
\end{equation*}
where the step-size $\tau$ satisfies the Armijo backtracking line-search along the projected directions \cite{nocedal2006numerical}.

\noindent Regarding the strong variation of material properties obtained in Section \ref{sec:OCP}, we add a regularizing weighting and force neighboring cells to have similar properties. The computed homogenized properties, indeed, refer to infinite repetition of equal unit cells, while in the most simple implementable configuration each hexagonal sub-domain is filled by a single unit cell which is thus surrounded by different ones. Limiting the difference between adjacent cells is thus beneficial for the equivalence of the behavior of the graded index metamaterial to the expected one.

\begin{figure}[t]
	\centering
	\includegraphics[width=\textwidth]{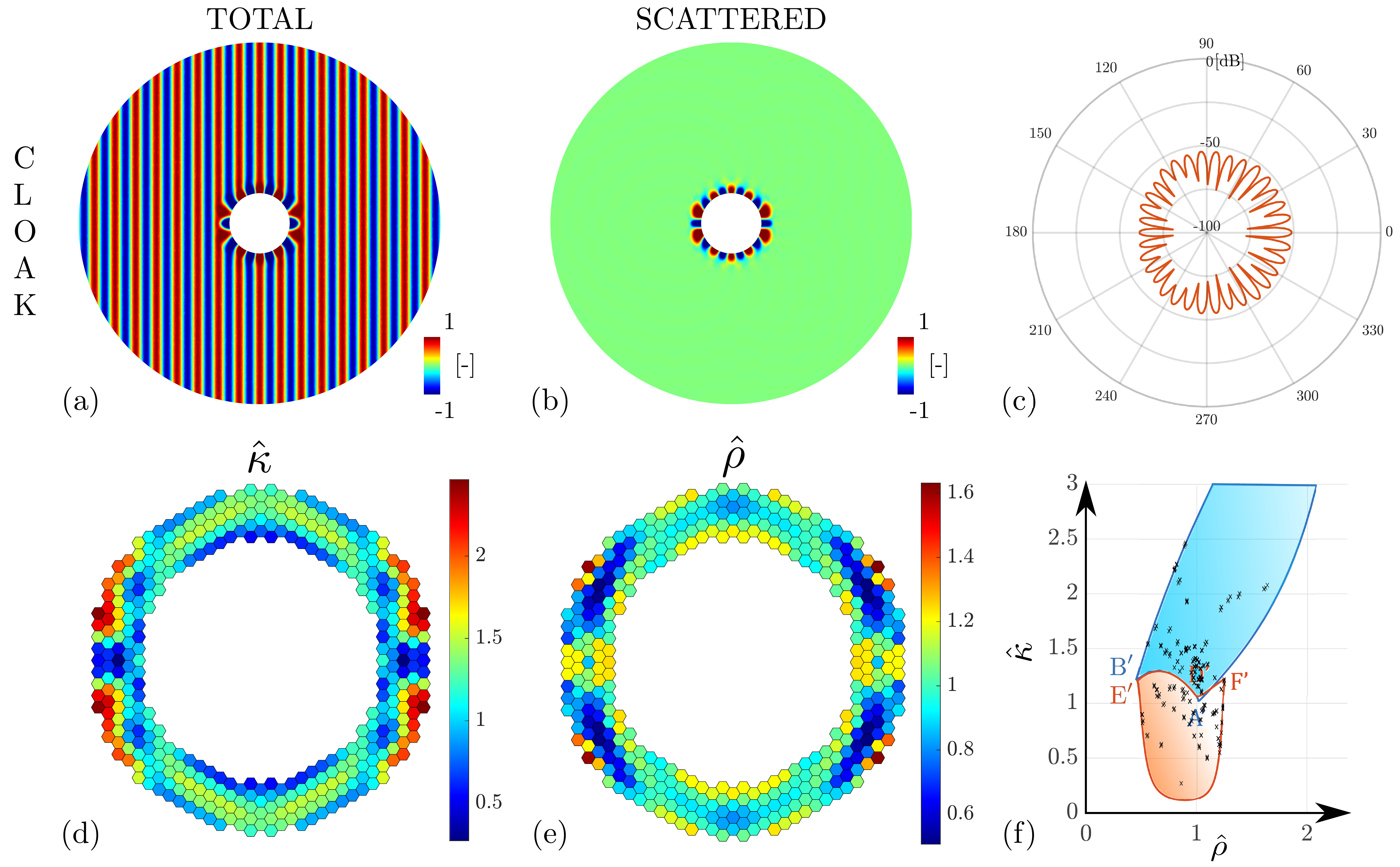}
	\caption[]{(a) Total pressure field computed with the cloak obtained with the material properties that are solution of the constrained OCP. (b) Associated scattered field. (c) Decibel reduction in scattered acoustic intensity with respect to the incident intensity. (d) Normalized bulk modulus distribution inside the cloak. (e) Normalized density distribution. (f) Each unit cell in the cloak represented as a $\hat{\rho} \times \hat{\kappa}$ pair falling inside the set $\mathcal{S}$.}
	\label{fig:cnstr scattering}
\end{figure}
\noindent Being the map from the geometrical parametrization to the equivalent properties regular and one-to-one, we can limit the geometrical dissimilarity between neighboring cells by constraining the equivalent properties on the $\hat{\rho} \times \hat{\kappa}$ plane or equivalently on the control space $\mathbf{u}$ and $\mathbf{v}$.

\noindent Thus, a penalty factor that weights the difference in the control intensity can be
 be introduced in the cost functional as : 
\begin{equation*}
	\frac{1}{2}\sum_{j=1}^{N_C}{\sum_{i\in \Lambda_j}{(u_i - u_j)^2}}
	 = \bb{u}^\top H \bb{u}
\end{equation*}
where $\Lambda_j$ is the set of cells adjacent to the $j$\textsuperscript{th} cell and $|\Lambda_j|$ its cardinality; the matrix $H$ is defined as:
\begin{equation*}
H_{ij}=
\begin{dcases}
    |\Lambda_i| & \text{if } i=j 
    \\
    -1 & \text{if } j\in \Lambda_i
    \\
    0 & \text{otherwise},
\end{dcases}
\end{equation*}
and we have used the identity $\sum_{j=1}^{N_c}{\sum_{i\in \Lambda_j}{u_i^2}} = \sum_{j=1}^{N_c}{|\Lambda_j|u_j^2}$. Note that the matrix $H$ corresponds to the Laplacian associated to the graph induced by the topology of the cells where an edge is present if the cells are neighbors. 
The graph is fully connected and thence its eigenvalues are nonnegative (see e.g.\ \cite{mesbahi2010graph}). The eigenvalue zero appears with multiplicity one and corresponds to the eigenvector space spanned by a vector of ones. Intuitively, this corresponds to the same control for all the cells.

As a result, the fully discrete cost function can be written as
\begin{equation*}
    J(\bb v,\bb u,\bb{p}_s) = \frac{\lambda_v}{2} \bb{v}^\top \left(H + D\right) \bb{v} + \frac{\lambda_u}{2} \bb{u}^\top \left(H + D\right) \bb{u} + \frac{1}{2}  \bb{p}^\dag M_{D_a}\bb{p},
\end{equation*}
where $D$ is the diagonal matrix whose entries are the areas of the associated cell. Due to the structure of $H$, it is clear that $H+D$ is positive definite. The fully discretized reduced gradients become:
\begin{equation*}
\begin{aligned}
\nabla J_{v_k} & = \lambda_v \left((H\bb v)_k + |D_{c,k}|v_k\,\right) +   \,e^{-v_k}  \bs{\lambda}^\dag\big( A_{k}\bb{p} + \bb{l}_k \big) 
\\[3pt]
\nabla J_{u_k} &=  \lambda_u \left((H \bb u)_k + |D_{c,k}|u_k\,\right)+e^{-u_k} \bs{\lambda}^\dag\big( B_{k}\bb{p} + \bb{d}_k \big).
\end{aligned}
\end{equation*}
The solution of the constrained optimization problem obtained by the PG method is shown in Figure~\ref{fig:cnstr scattering}. In particular, in Figure~\ref{fig:cnstr scattering}(a)-(b) and (c) are depicted the total field, the scattered field, and the polar dependence of the decibel gain in scattered intensity computed with respect of the incident intensity, as previously done in the unconstrained scenario. 

The performance in terms of scattering reduction are comparable to those obtained without the constraints. Moreover, Figure~\ref{fig:cnstr scattering}(d)-(e) show the obtained solution of the constrained optimization in terms of material properties distribution, i.e.\ the normalized bulk modulus and density, respectively. Finally, Figure \ref{fig:cnstr scattering}(f) show the location of each unit cell as black markers in the $\hat{\rho}\times \hat{\kappa}$ plane. Note that the obtained material properties lie inside the reachable set $\mathcal{S}$ or on its boundary $\partial\mathcal{S}$ whenever the feasibility constraint is active.

\section{Design of the Microstructured Cloak and Validation}

\begin{figure}[t]
	\centering
	\includegraphics[width=\textwidth]{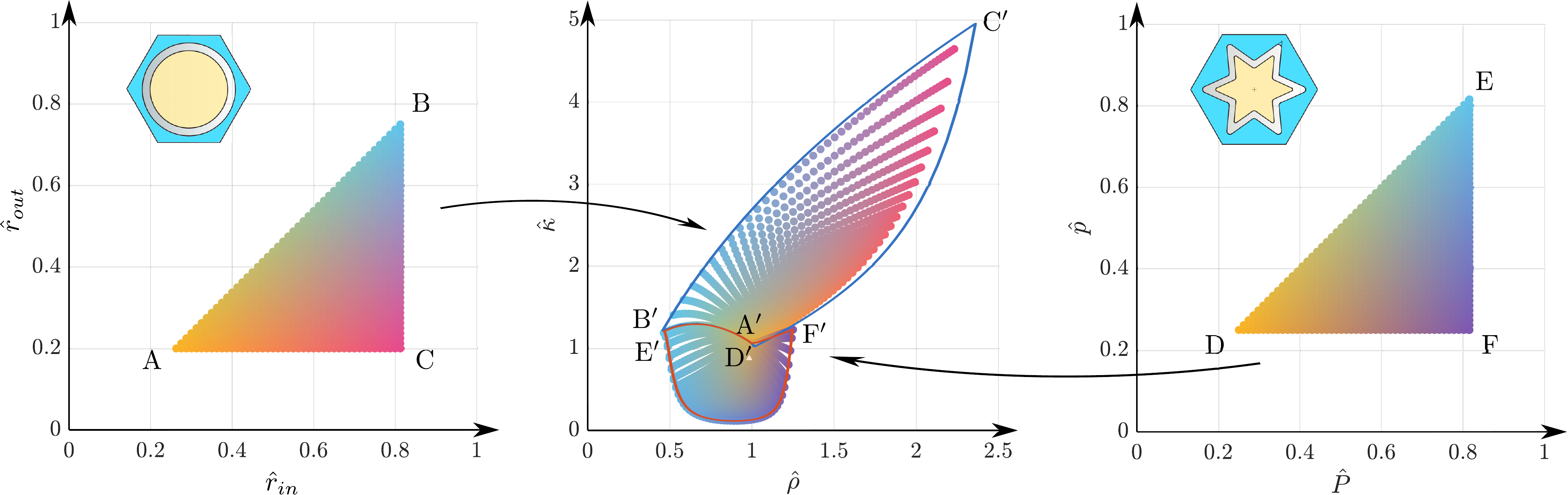}
	\caption[]{Direct mapping between the space of geometrical parameters and the space of homogenized material properties, that is inverted to solve for the design of the microstructure once the constrained optimal control problem is solved. Colors help to trace visually each different point from the space of homogenized material properties back to that of its geometrical features.  }
	\label{fig:cloak1}
\end{figure}

Once the optimal required material properties are found, the inverse engineering problem of finding the microstructure geometry that exhibit those $\hat{\rho}$ and $\hat{\kappa}$ pairs has to be solved. This being a much more difficult problem than the direct one, it is usually tackled adopting optimization algorithms, either parametric of evolutionary, that employ as cost function the distance between the required desired material properties and those obtained by homogenization on the considered lattice \cite{chen2017broadband, quadrelli2021experimental}. In the case at hand, the simplicity of the geometry of the considered unit cells, which is univocally determined in both configurations by a pair of parameters, allows for a direct mapping of the whole $\mathcal{C}^\odot$ and $\mathcal{C}^\star$ spaces into the $\hat{\rho} \times \hat{\kappa}$ one. Once this map is computed, it can subsequently be used to solve the inverse engineering; in particular, the homogenized material properties are computed for the grid of points shown on the $\hat{r}_{in} \times \hat{r}_{out}$ and $\hat{p}\times \hat{P}$ spaces in Figure~\ref{fig:cloak1} and the resulting discrete map is used for a first guess of the cells geometrical parameters when required $(\hat{\rho}, \hat{\kappa})$ values are specified. An optimization routine allows then to refine the properties of each cell with few iterations.

\begin{figure}[t]
	\centering
    \subfloat[][]
    {\includegraphics[width=.48\textwidth]{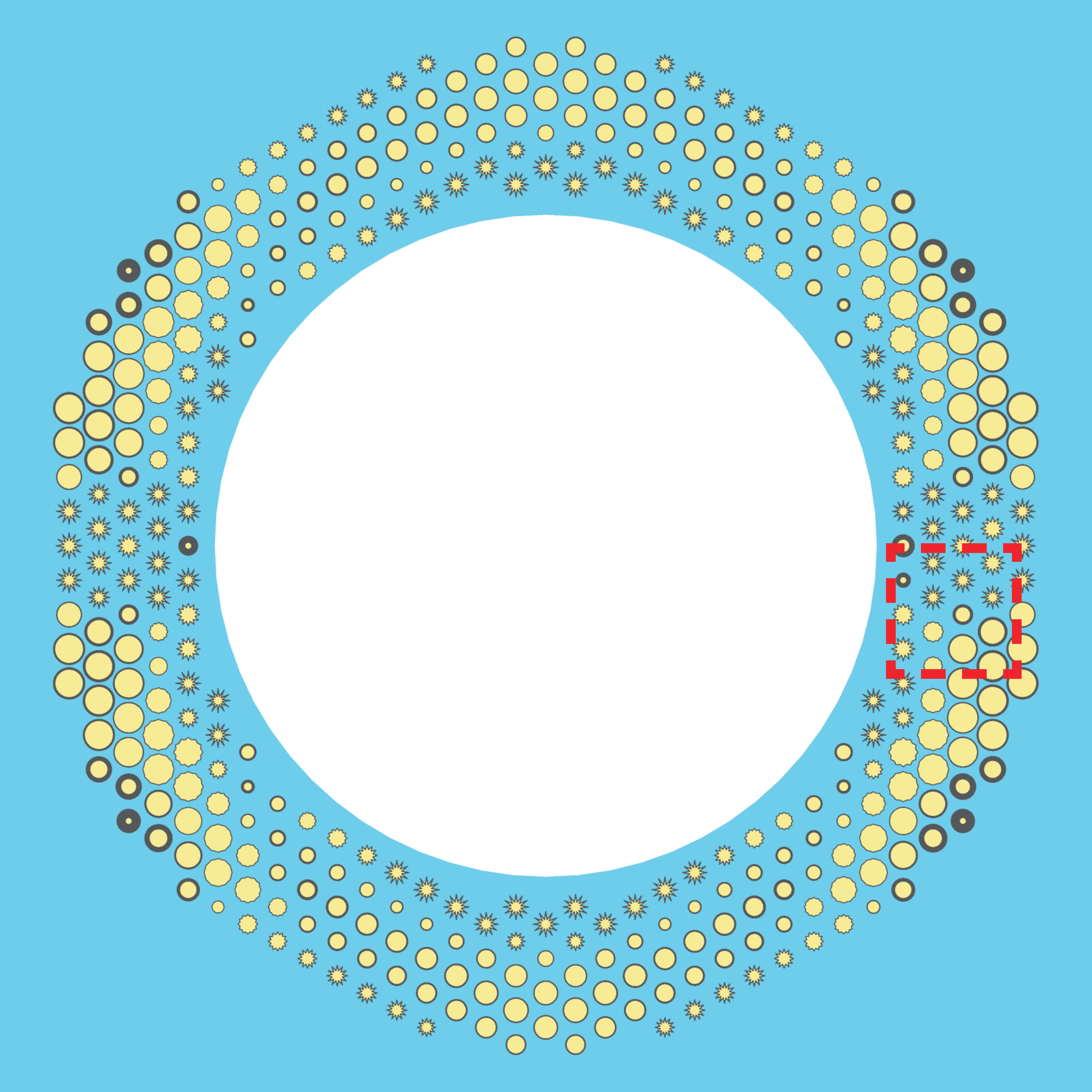}}
	\quad
    \subfloat[][]
    {\includegraphics[width=.48\textwidth]{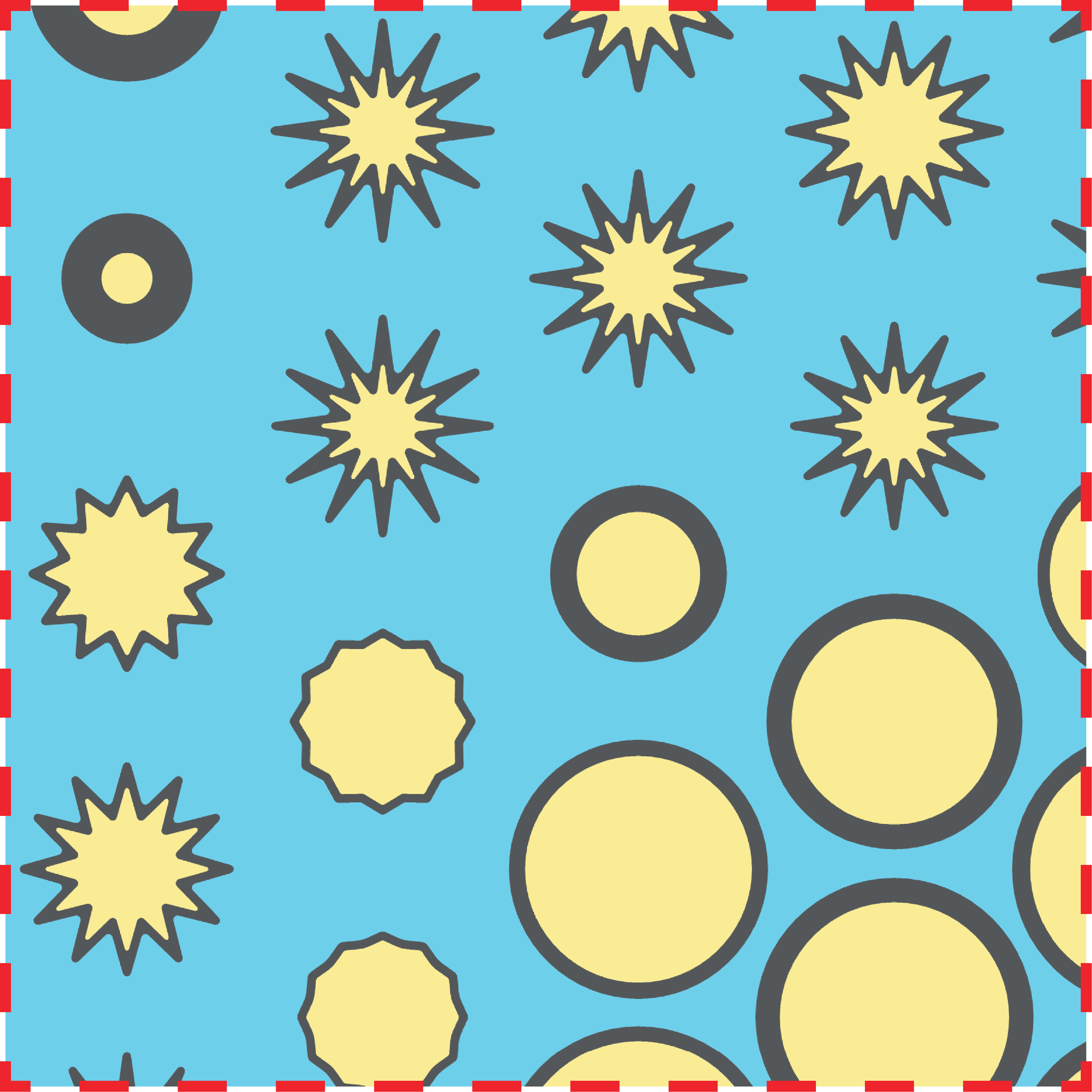}}
    \caption[]{(a)~Schematic of the geometry of the cloak made by the hexagonal lattice of inclusions. Aluminium part are depicted in grey, air in yellow and water in light blue. (b)~Magnification of the red box depicted on (a) that shows the geometry of the inclusions.}
	\label{fig:geom circ}
\end{figure}

Following the aformentioned design procedure, the entire cloak geometry obtained from the solution of the constrained OCP is defined and the resulting microstructure is depicted in Figure~\ref{fig:geom circ}, where colors are used to distinguish between domains filled by air, aluminium or water. A fully coupled structural/acoustic frequency domain finite element simulation of the designed cloak is carried out by means of the commercial software \Comsol{}, in order to test its performances when considering the actual implemented structure. 

\begin{figure}[t]
	\centering
	\includegraphics[width=\textwidth]{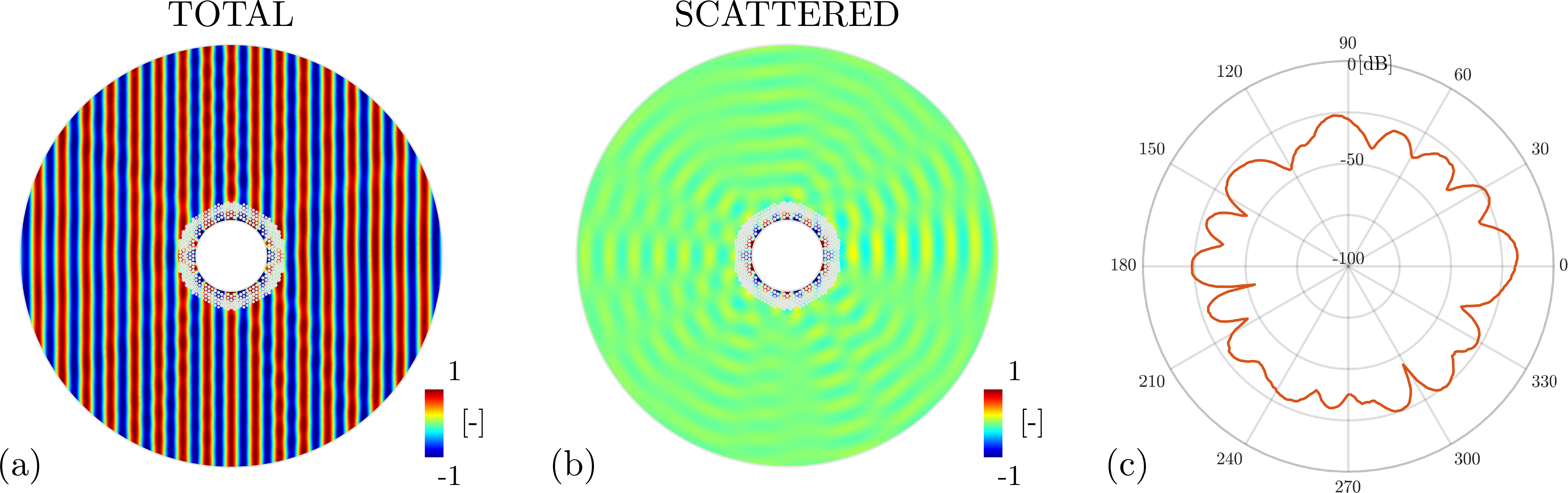}
	\caption[]{Fully coupled structural acoustic finite element simulation of the microstructured cloak. (a)~Total pressure field (b)~Associated scattered pressure field. (c)~Decibel reduction in scattered acoustic intensity.  }
	\label{fig:cloak2}
\end{figure}

A first-order approximation of the Sommerfeld absorbing condition is assigned on the boundary $\Gamma_e$ to approximate an unbounded domain.
The results are shown in Figure~\ref{fig:cloak2}.

By looking at the scattered intensity plot of Figure~\ref{fig:cloak2}(c), we can state that the outgoing energy is two orders of magnitude lower with respect to the uncloaked case; then the obstacle is undetectable. The discrepancies with respect to the simulation performed with the homogenized properties can be attributed to the fact that one single unit cell has been considered to fill each cloak sub-domain, while an infinite microstructure should ideally be placed there instead.

\begin{figure}[t]
	\centering
	\includegraphics[width=\textwidth]{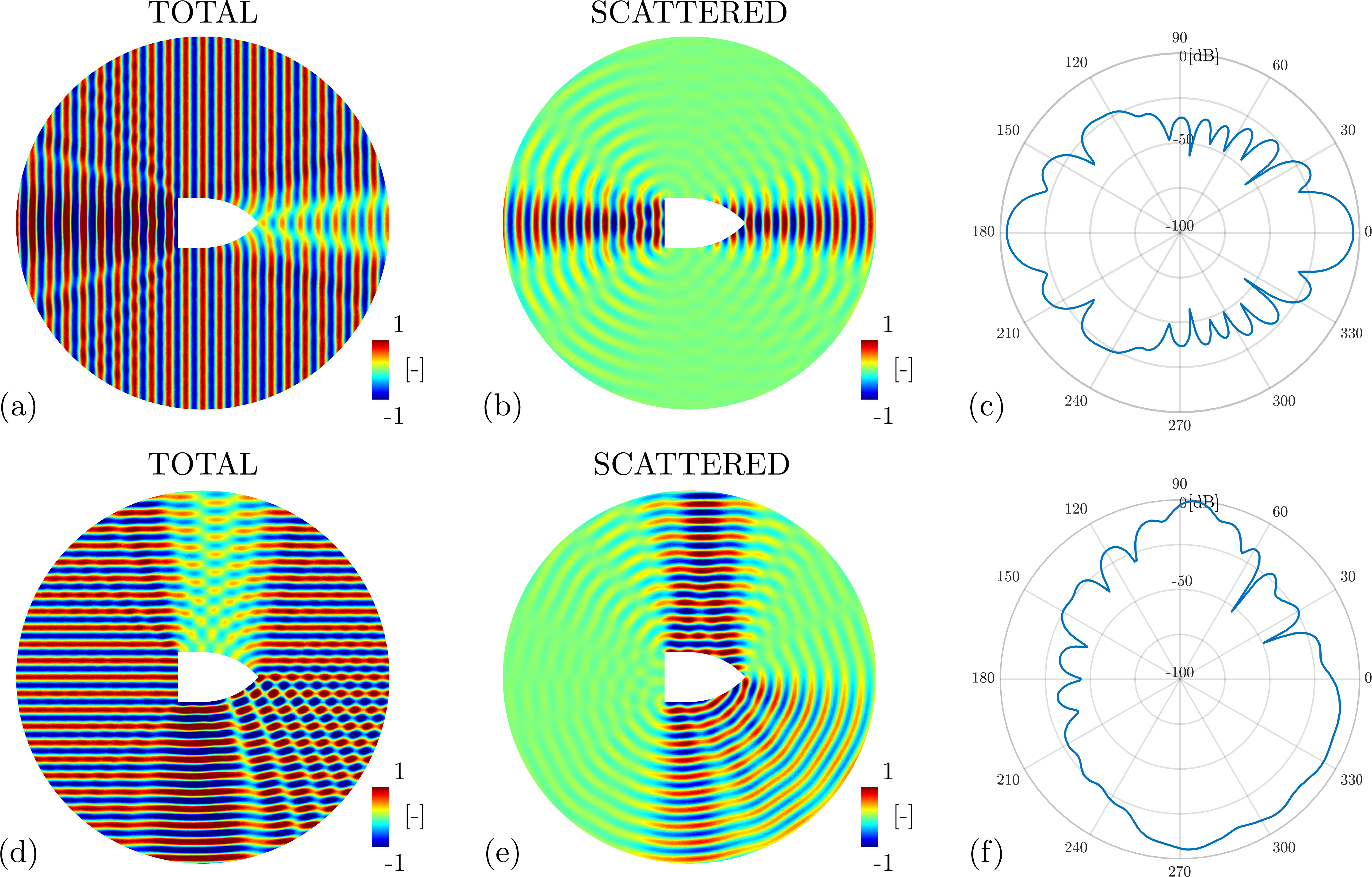}
	\caption[]{(a) Total pressure field for a wave incident from left with wavelenght $\lambda_1=\SI{20.4}{\percent}$ of the ship characteristic length. (b)~Associated scattered pressure field. (c)~Decibel reduction in scattered acoustic intensity with respect to the incident intensity. (d)~Total pressure field for a wave incident from bottom with wavelenght $\lambda_2=\SI{18.9}{\percent}$ of the ship characteristic length. (e)~Associated scattered pressure field. (f)~Decibel reduction in scattered acoustic intensity with respect to the incident intensity.}
	\label{fig:cloak3}
\end{figure}

As a further test case, a constrained OCP is set to find the optimal material properties' distribution to cloak the silhouette of a ship, i.e.\ an obstacle with a non axisymmetric contour. The probing acoustic field consists in the superposition of an incident plane wave with wavelength $\lambda_1=\SI{20.4}{\percent}$ of the ship characteristic length $L$ and direction $\mathbf{a}=[1,0]$ (horizontal incidence) and a plane wave with wavelength $\lambda_2=\SI{18.9}{\percent}\,L$ and direction $\mathbf{a}=[0,1]$ (vertical incidence).  The size of each hexagonal sub-domain is $\SI{8.7}{\percent}\,\lambda_1=\SI{9.3}{\percent}\,\lambda_2$.  Figure~\ref{fig:cloak3} shows the uncloaked case scenario in terms of total fields, scattered fields and scattered intensity for both horizontal and vertical incidence. This choice for the probing incident field allows to test the performance of the method when multiple frequencies and directions are taken into account: for this reason, the definition of the OCP is modified as follows. Let us consider a number $N_f$ of incident pressure fields $p_{i,h}$, $h\in\{1,\dots,N_f\}$. The governing equations are linear with respect to the pressure, thus the superposition principle holds and we can modify the objective functional by weighting the sum of the scattered fields for each probing frequency. Indeed, we can select as $J$:
\begin{equation*}
J_{N_f} = 
    {\lambda_v \over 2}\int_{D_c}{v^2\,d\Omega} + {\lambda_u\over 2} \int_{D_c}{u^2\,d\Omega} + 
	\frac{1}{2}\sum_{h=1}^{N_f}{ \int_{D_a}{\bar{p}_{s,h} p_{s,h}	\,d\Omega}}
\end{equation*}
where each scattered pressure $p_{s,h}$ satisfies the state dynamics \eqref{eq:state_dyn} with frequency $\omega_h$ and forcing terms determined by $p_{i,h}$. Note that the PDE constraints are now $N_f$. With similar arguments as for the previous section, we can form a Lagrangian which comprises the sum of the PDE constraints. From the latter, we can compute $N_f$ adjoint equations of the form \eqref{eq:PDE adjoint} where the right-hand side depends on $p_{s,h}$ only. Note that the control functions are the same for each state and adjoint equation. In this way, the reduced gradients can be computed as:

\begin{equation*}
\begin{aligned}
    \nabla J_{\bb{v}} &= \lambda_v \Big(\int_{D_c} \bs{\psi} \bs{\psi}^{\top} \, d\Omega\Big) \bb{v} +  \sum_{h=1}^{N_f}\Re\Big\{\int_{D_c} \bs{\psi} \, a \nabla(p_{s,h} + p_{i,h})\cdot\nabla \bar{\lambda}_h \, d\Omega\Big\}
    \\
    \nabla J_{\bb{u}} &= \lambda_u \Big(\int_{D_c} \bs{\psi} \bs{\psi}^{\top} \, d\Omega\Big) \bb{u} + \sum_{h=1}^{N_f}\Re\Big\{\int_{D_c} \bs{\psi} \, b\omega_h^2 (p_{s,h} + p_{i,h}) \bar\lambda_h \, d\Omega \Big\}.\\
\end{aligned}
\end{equation*}

\begin{figure}[t]
	\centering
	\includegraphics[width=\textwidth]{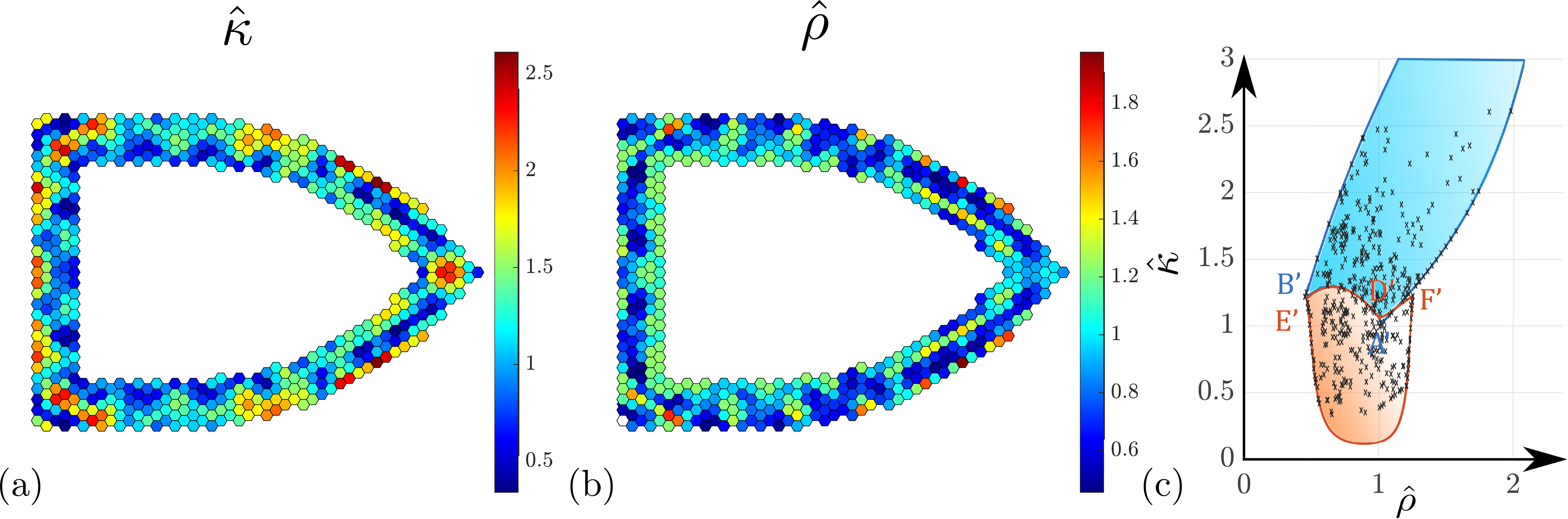}
	\caption[]{(a)~Normalized bulk modulus distribution inside the cloak. (b)~Normalized density distribution. (c)~Each unit cell in the cloak represented as a $\hat\rho\times\hat\kappa$  pair falling inside the set $\mathcal S$}
	\label{fig:cloak4}
\end{figure}

\begin{figure}[t]
	\centering
	\includegraphics[width=\textwidth]{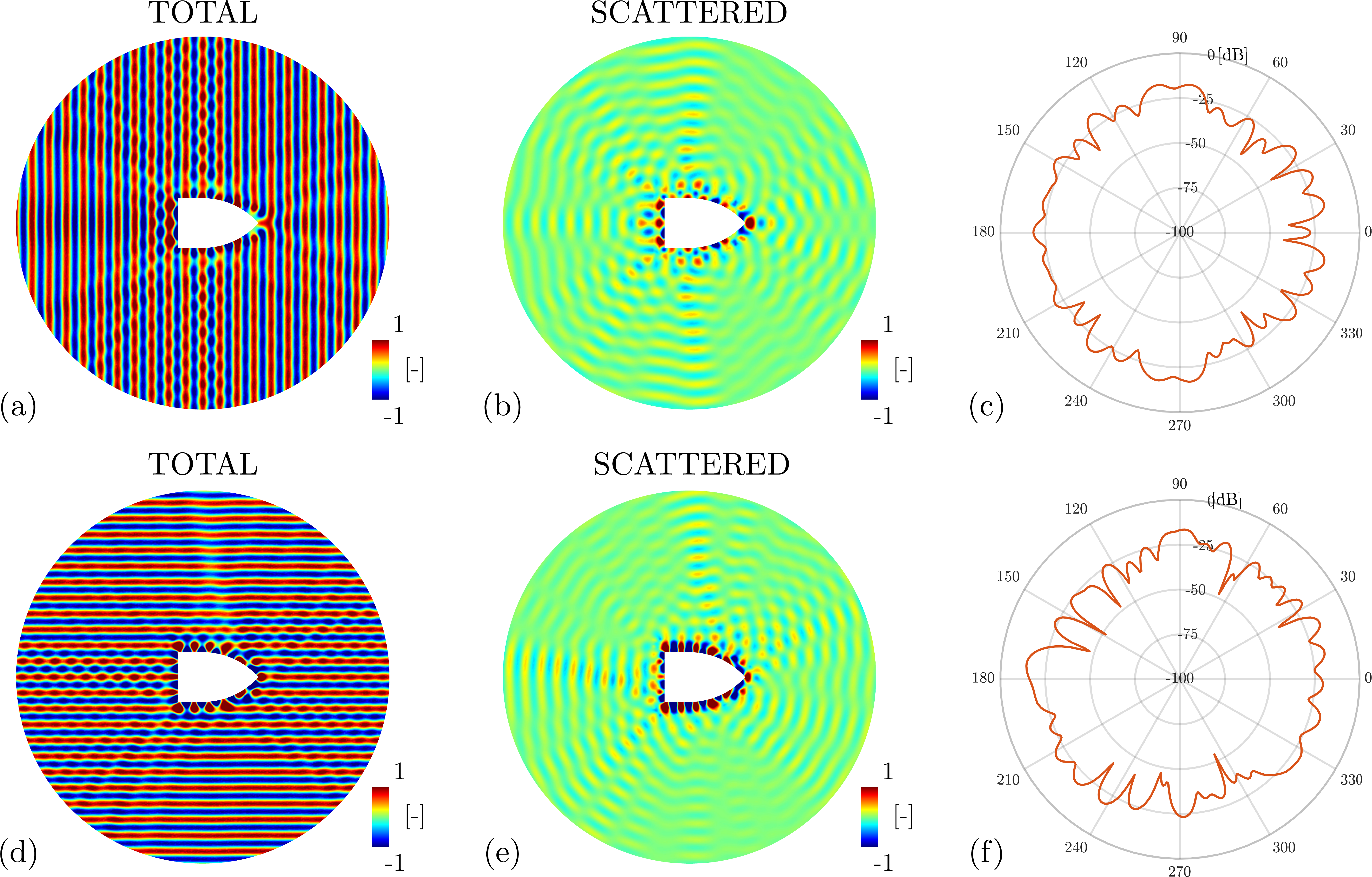}
	\caption[]{Acoustic fields obtained with the cloak made by the properties depicted in Figure~\ref{fig:cloak4} (a)~Total pressure field for horizontal incidence (b)~Associated scattered pressure field. (c)~Decibel reduction in scattered acoustic intensity. (d)~Total pressure field for vertical incidence. (e)~Associated scattered pressure field. (f)~Decibel reduction in scattered acoustic intensity.}
	\label{fig:cloak5}
\end{figure}
\begin{figure}
\centering
\subfloat[][]
{\includegraphics[width=.4\textwidth,trim=0 0 20 0]{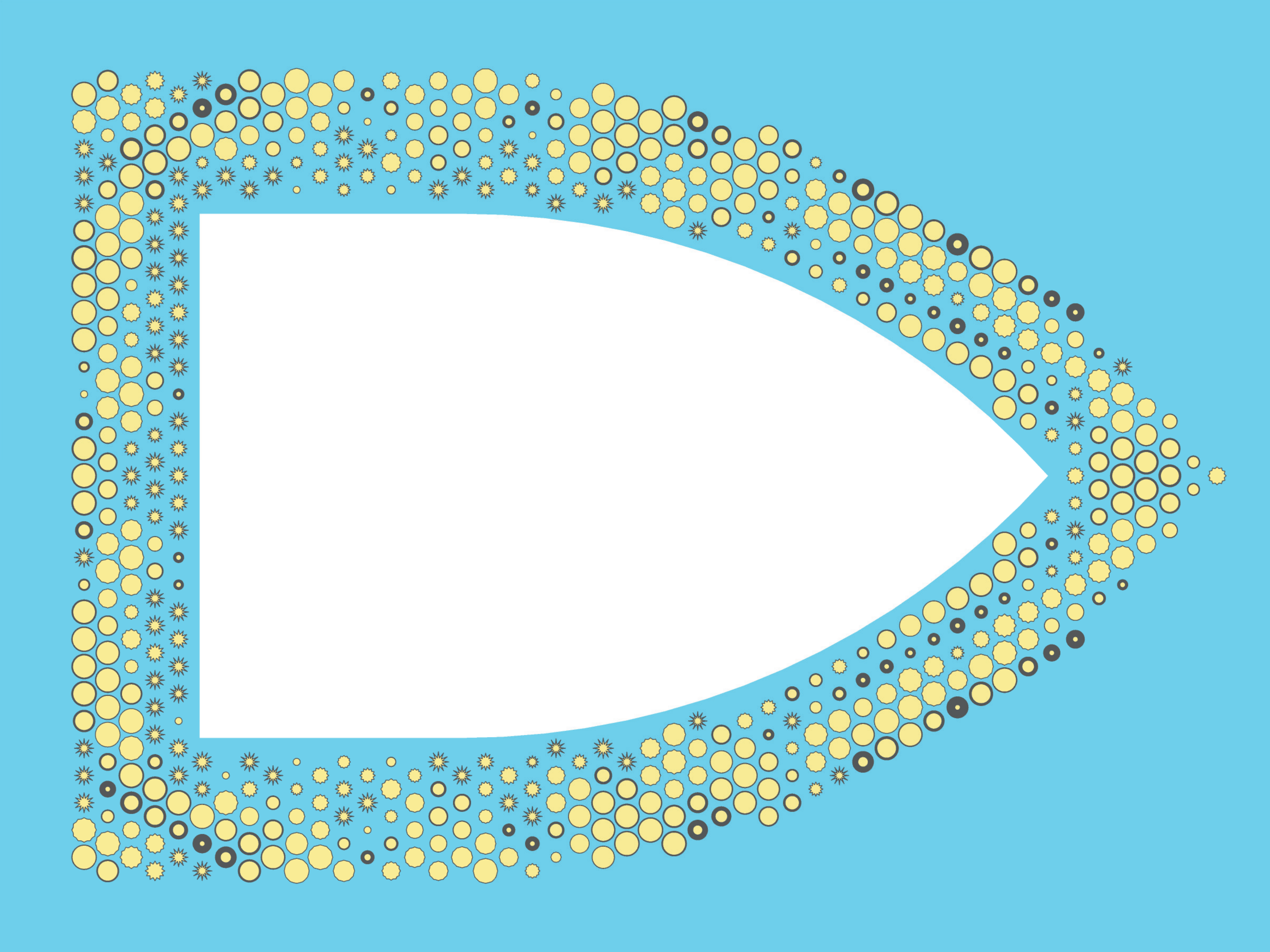}}
\quad
\subfloat[][]
{\includegraphics[width=.55\textwidth,trim=20 20 20 0]{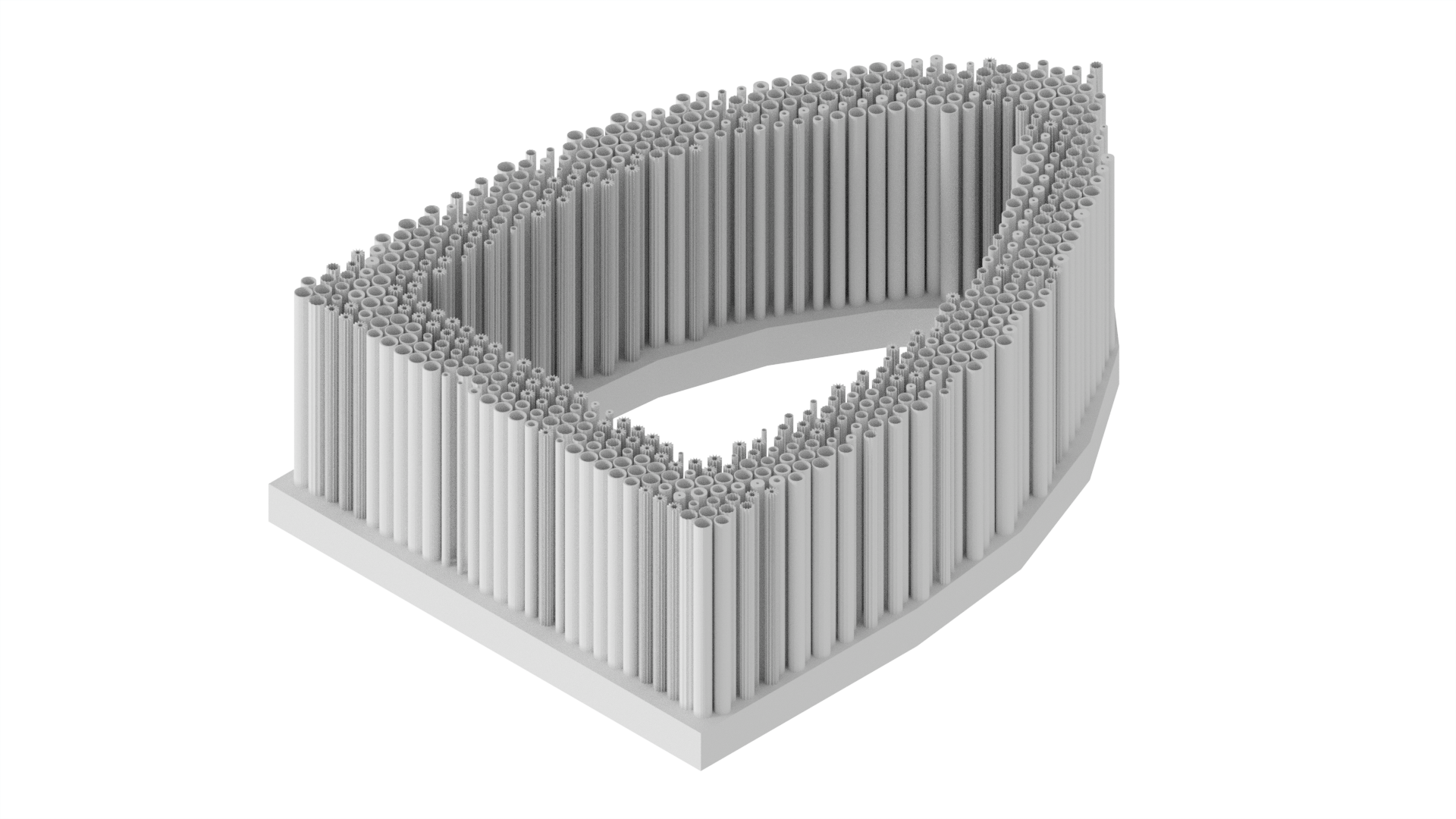}}
\caption{(a)~Schematic of the entire cloak microstructure geometry for the ship obstacle. Grey is used for the Aluminium inclusions, yellow for air and light blue for water. (b)~3D render of the extruded geometry useful for experimental validation.}
\label{fig:geom barca}
\end{figure}

\begin{figure}[t]
	\centering
	\includegraphics[width=\textwidth]{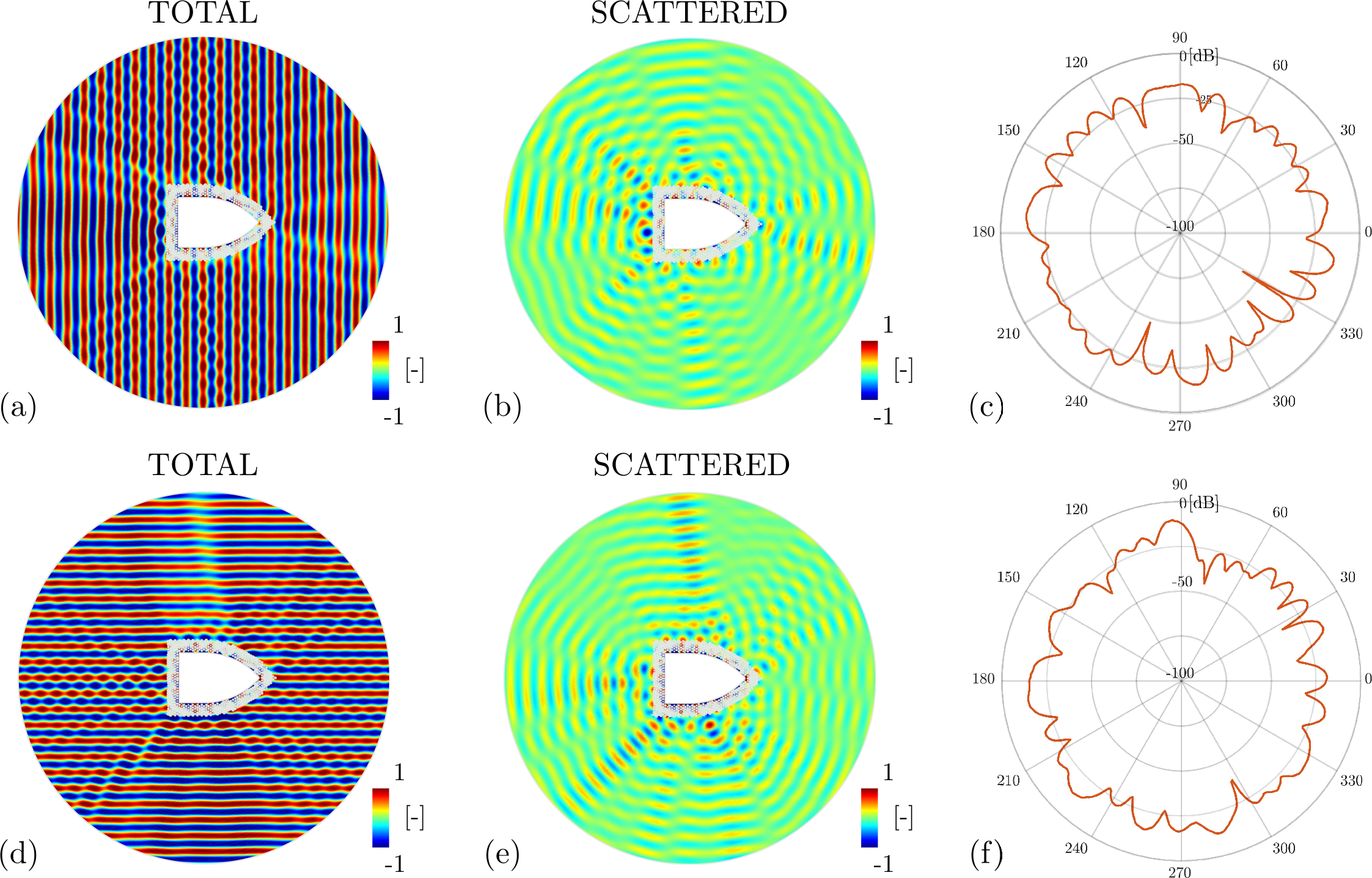}
	\caption[]{Fully coupled structural acoustic finite element simulation of the ship equipped with microstructured cloak. (a)~Total pressure field for horizontal incidence (b)~Associated scattered pressure field. (c)~Decibel reduction in scattered acoustic intensity. (d)~Total pressure field for vertical incidence. (e)~Associated scattered pressure field. (f)~Decibel reduction in scattered acoustic intensity.}
	\label{fig:cloak6}
\end{figure}

The multi-frequency problem is solved with the PG method and the results are shown in Figure~\ref{fig:cloak4}(a,b) in terms of material properties distributions while in Figure~\ref{fig:cloak4}(c) it is shown that they all lie in the feasible set $\mathcal{S}$. The corresponding acoustic fields and scattered intensity are shown in Figure~\ref{fig:cloak5} for comparison with the uncloaked scenario. A \SI{25}{[\decibel]} reduction of scattered intensity is obtained in both the backward and forward scattering directions with respect to the uncloaked case. As done for the axisymmetric case, the geometry of the actual cloak is deduced optimizing for each single unit cell, and the resulting design is depicted in Figure~\ref{fig:geom barca}(a). Such geometry can be easily extruded in the third dimension and fabricated with state of the art 3D printing technologies for experimental validations (Figure~\ref{fig:geom barca}(b)). In order to provide an accurate numerical validation of the cloak design, a fully coupled acoustic/structure simulation is performed with the commercial software \Comsol{}.

Figure~\ref{fig:cloak6} shows the computed fields and figures of merit of the ship equipped with the cloak composed by the actual microstructure, showing good agreement between the obtained performances and the ideal ones obtained with the homogenized material properties (Figure~\ref{fig:cloak5}).

\section{Conclusions}

In this paper, we have introduced a general acoustic cloaking design strategy that simultaneously aims at reducing the complexity of the required microstructures and enlarge the set of geometries that can be cloaked with respect to traditional Transformation-based methods. This is achieved by synergic use of PDE-constrained optimization, to find the isotropic material distribution that minimizes scattering, and parametric structural optimization, to design simple hexagonal lattices of inclusions that match the required densities and bulk moduli. More than that, such two scale-optimization problem is formulated in such a way that the two stages, i.e. the computation of the \textit{macroscale} material properties distribution and the \textit{microscale} design, are not disconnected steps but intimately linked together, in order to retain the optimality of the solution found. This is done at the OCP level by considering as control space a suitable linear combination of indicator functions which corresponds to the topology of the hexagonal lattices, and constraining the controls to take values inside a feasible region that is pre-computed analyzing all the possible considered unit cell geometries. The method is tested against the usual axisymmetric cloaking scenario, producing a two orders of magnitude mean reduction of intensity over the whole azimuthal scattering directions. Then, a more complicated scenario is considered, where an arbitrary shaped obstacle is probed by two acoustic waves with different frequency and incident direction. The solution is found to reduce the backscattered and forward scattered wave with performance comparable to those obtained in the simple axisymmetric scenario. With the simplicity of the considered geometries, this manuscript paves the way for experimental validation of the acoustic cloaking principle with arbitrary obstacle shapes.

\bibliographystyle{RS}
\bibliography{references.bib}

\end{document}